\documentclass[a4paper,11pt]{amsart}
\usepackage{amssymb}
\usepackage{amscd}

\newtheorem{theorem}{Theorem}[section]
\newtheorem{lemma}[theorem]{Lemma}

\newtheorem{proposition}[theorem]{Proposition}
\newtheorem{corollary}[theorem]{Corollary}

\theoremstyle{definition}

\theoremstyle{remark}

\usepackage{graphicx}
\usepackage{amsmath}
\usepackage{amsfonts}
\usepackage{amssymb}
\newcommand{\be}{\begin{equation}}

\newcommand{\ee}{\end{equation}}

\newcommand{\cB}{{\mathcal B}}

\newcommand{\cG}{\mbox{$\mathcal{G}$}}

\newcommand{\si}{\sigma}

\newcommand{\ba}{\begin{array}}

\newcommand{\ea}{\end{array}}

\newcommand{\beq}{\begin{eqnarray}}

\newcommand{\eeq}{\end{eqnarray}}

\newtheorem{lm}{lemma}

\newtheorem{thee}{theorem}

\newtheorem{proo}{proposition}

\newtheorem{co}{corollary}

\newtheorem{rem}{remark}

\newtheorem{deff}{definition}

\newcommand{\bd}{\begin{deff}}

\newcommand{\ed}{\end{deff}}

\newcommand{\bl}{\begin{lm}}

\newcommand{\el}{\end{lm}}

\newcommand{\bp}{\begin{proo}}

\newcommand{\ep}{\end{proo}}

\newcommand{\bt}{\begin{thee}}

\newcommand{\et}{\end{thee}}

\newcommand{\bc}{\begin{co}}

\newcommand{\ec}{\end{co}}

\newcommand{\brm}{\begin{rem}}

\newcommand{\erm}{\end{rem}}

\usepackage{t1enc}
\def\frak{\mathfrak}
\def\Bbb{\mathbb}
\def\Cal{\mathcal}

\newcommand{\newc}{\newcommand}

\renewcommand{\exp}{\operatorname{exp}}

\renewcommand{\o}{\circ}

\let\ccdot\cdot
\def\cdot{\hbox to 2.5pt{\hss$\ccdot$\hss}}

\newc{\aR}{\mbox{\boldmath{$ R$}}}
\newc{\aS}{\mbox{\boldmath{$ S$}}}
\newc{\aT}{\mbox{\boldmath{$ T$}}}
\newc{\aW}{\mbox{\boldmath{$ W$}}}

\newcommand{\aI}{\mbox{\boldmath{$ I$}}}

\newcommand{\aH}{\mbox{\boldmath{$ H$}}}

\newc{\aK}{\mbox{\boldmath{$ K$}}}
\newc{\aL}{\mbox{\boldmath{$ L$}}}

\newcommand{\ce}{{\Cal E}}

\usepackage{amssymb}
\usepackage{amscd}

\newcommand{\nd}{\nabla}
\renewcommand{\and}{\boldsymbol{\nd}}

\newcommand{\ul}[1]{\underline{#1}}

\newcommand{\cT}{{\mathcal T}}

\let\i=\iota

\newcommand{\aM}{\boldsymbol{M}}

\newcommand{\cW}{{\Cal W}}

\newcommand{\nn}[1]{(\ref{#1})}

\newc{\obstrn}[2]{B^{#1}_{#2}}

\newcommand{\rpl}                         
{\mbox{$
\begin{picture}(12.7,8)(-.5,-1)
\put(0,0.2){$+$}
\put(4.2,2.8){\oval(8,8)[r]}
\end{picture}$}}

\newcommand{\lpl}                         
{\mbox{$
\begin{picture}(12.7,8)(-.5,-1)
\put(2,0.2){$+$}
\put(6.2,2.8){\oval(8,8)[l]}
\end{picture}$}}

\usepackage{ifthen}

\newc{\tensor}[1]{#1}
\newc{\Mvariable}[1]{\mbox{#1}}
\newc{\down}[1]{{}_{#1}}
\newc{\up}[1]{{}^{#1}}

\newc{\JulyStrut}{\rule{0mm}{6mm}}
\newc{\midtenPan}{\mbox{\sf S}}
\newc{\midten}{\mbox{\sf T}}
\newc{\midtenEi}{\mbox{\sf U}}
\newc{\ATen}{\mbox{\sf E}}
\newc{\BTen}{\mbox{\sf F}}
\newc{\CTen}{\mbox{\sf G}}

\def\sideremark#1{\ifvmode\leavevmode\fi\vadjust{\vbox to0pt{\vss
 \hbox to 0pt{\hskip\hsize\hskip1em
 \vbox{\hsize3cm\tiny\raggedright\pretolerance10000
 \noindent #1\hfill}\hss}\vbox to8pt{\vfil}\vss}}}%
                        
\numberwithin{equation}{section}

\begin{document}
\renewcommand{\today}{}
\title{Projective BGG equations, algebraic sets, and compactifications of Einstein geometries}

\author{A.\ \v Cap, A.R.\ Gover, M.\ Hammerl}

\address{A.\v C. and M. H.: Faculty of Mathematics\\
University of Vienna\\
Nordbergstr. 15\\
1090 Wien\\
Austria\\
A.R.G.:Department of Mathematics\\
  The University of Auckland\\
  Private Bag 92019\\
  Auckland 1142\\
  New Zealand;\\
Mathematical Sciences Institute\\
Australian National University \\ ACT 0200, Australia} \email{Andreas.Cap@univie.ac.at}
\email{r.gover@auckland.ac.nz}
\email{Matthias.Hammerl@univie.ac.at}

\vspace{10pt}

\renewcommand{\arraystretch}{1}


\begin{abstract} 
  For curved projective manifolds we introduce a notion of a normal
  tractor frame field, based around any point. This leads to canonical
  systems of (redundant) coordinates that generalise the usual
  homogeneous coordinates on projective space.  These give preferred
  local maps to the model projective space that encode geometric
  contact with the model to a level that is optimal, in a suitable sense.  In
  terms of the trivialisations arising from the special frames, normal
  solutions of classes of natural linear PDE (so-called first BGG
  equations) are shown to be necessarily polynomial in the generalised
  homogeneous coordinates; the polynomial system is the pull back of  
  a  polynomial system that solves the corresponding problem on
  the model. Thus questions concerning the zero locus of solutions, as
  well as related finer geometric and smooth data, are reduced to a
  study of the corresponding polynomial systems and algebraic sets. We
  show that a normal solution determines a canonical manifold
  stratification that reflects an orbit decomposition of the model.
  Applications include the construction of structures that are
  analogues of Poincar\'e-Einstein manifolds.
\end{abstract}

\maketitle
\renewcommand{\arraystretch}{1.5}

\pagestyle{myheadings} \markboth{\v Cap, Gover, Hammerl}{Projective
  equations, algebraic sets, and compactifications}

\thanks{ARG gratefully acknowledges support from the Royal Society of
  New Zealand via Marsden Grant no.\ 06-UOA-029; A\v C and MH
  gratefully acknowledge support by project P19500--N13 of the ``Fonds
  zur F\"orderung der wissenschaftlichen Forschung'' (FWF) and the
  hospitality of the University of Auckland.}

\maketitle
\section{Introduction}

Natural (or geometric) partial differential equations (PDE) are those
which, in a suitable sense, are determined by some underlying
geometry.  Given such an equation, an important problem is to expose
the geometric content of a solution, and the implications of its
existence.  Example questions close to our focus here are the
following.  What general results can be established concerning the
nature and geometric structure of the solution's zero locus?  What is
the relationship of these features to the ambient structure? At the
most primitive level there are questions of topology and smoothness;
at the next level, it can be that the zero locus satisfies an
interesting smooth embedding equation and inherits a rich intrinsic
geometric structure.

Here we study a large class of solutions to overdetermined PDE arising
naturally in projective geometry.  We show that, to a surprising
extent, these problems can be reduced to vastly simpler questions of
an algebraic geometric type. This leads to a conceptual and practical
way to describe, and manage geometrically, compactifications of curved
geometries; the compactifications involved are naturally related to the
geodesic structure.

Our work is partly inspired by the rich programme surrounding
Poincar\'e-Einstein (PE) manifolds; these were introduced by
Fefferman-Graham as a tool for constructing conformal invariants
\cite{FGast}. A PE structure consists of a conformal compactification
of a geodesically complete Einstein-pseudo-Riemannian manifold, and
broadly the programme involves naturally relating geometry and fields
on a conformal boundary with Einstein-Riemannian geometry and field
theory on the interior. Over the past decade this problem has been a
domain of intense interest and deep progress
\cite{AndReg,CQY,FGrQ,FGnew,GrZ,lee,MP}, in part the work has been driven by
strong links with the AdS/CFT correspondence of Maldacena
\cite{Maldacena,GKP}, see e.g.\ \cite{GrSrni,GrW,SS}. The original PE
construction was preceded by a K\"ahler-Einstein-CR geometry analogue
\cite{CY,Feff} and recently there has been work to extend the picture
to quaternionic-K\"ahler metrics and beyond \cite{Biquard,BiqMazz}.

It was observed in \cite{GoPrague} that a PE manifold is  the same as a
conformal manifold equipped with a solution of a certain conformally
invariant PDE; this solution having the property that its (necessarily
smooth and embedded) zero locus is precisely the topological
boundary. That result led to an effective approach to certain key
problems for these structures, extension to the notion of almost
Einstein manifolds \cite{GoMGC,Go-al}, and also methods for
geometrically constructing, and partly characterising, examples of PE
manifolds \cite{GoLeitprog}. In \cite{Go-al} it is  seen that the
almost Einstein class also naturally includes asymptotically locally
Euclidean (ALE) structures that admit isolated point conformal
compactification; in fact the nature of the compactification is shown
to be an easy consequence of the compatibility of Ricci-flatness with the
governing conformal PDE.

Here we show that considering a class of solutions to similar, but
essentially different equations, leads to a natural extension of these
ideas.  In fact the scope is broader than this suggests, as the wider
perspective draws in considerable new phenomena.  In the simplest
class of cases the linear equations studied take the form
\begin{equation}\label{scales}
  \big(\nabla_{(a_1}\nabla_{a_2}\cdots \nabla_{a_{k+1})}  + 
\mbox{lower order terms}\big)\si =0,
\end{equation}
where $\nabla$ is an affine connection, $\si$ is a function and the
$(\cdots)$ indicates taking the symmetric part over the enclosed
indices. (Here, and in many places throughout, we use Penrose's
abstract index notation \cite{ot}.) These equations have a symmetry
known as projective invariance that plays an important role.  We see
in Section \ref{prelim} that certain solutions for the cases $k=1$ and
$k=2$ do indeed lead to structures that are analogous to PE
manifolds. In particular in Section \ref{KleinE} the $k=2$ case yields
a curved analogue of the hyperbolic ball that has been earlier described from
quite a different perspective as a case of a ``projectively compact
metric'' in \cite{FGnew}. 
This is a manifold with boundary. As in the
case of a PE manifold, the boundary has a canonical conformal
structure and the interior has a geodesically complete (Riemannian)
negative Einstein structure. However in this case a projective
compactification is involved, which emphasises the role of geodesics;
this is strictly different from conformal compactification thus the
structure is {\em not} PE, see Proposition \ref{KEP}.  There are
analogues for all signatures and the Lorentzian case should be of
interest to the general relativity community.  These examples also
show rather clearly that although the equations we consider are
linear, the integrability conditions for these can involve very
interesting non-linear conditions (such as the Einstein equations), on
which we obtain a new perspective.

Let us now be specific about the full class of equations we treat.  We
restrict our attention to natural equations on a projective manifold
of dimension at least 2. Recall that this consists of a manifold $M$
equipped with an equivalence class $p$ of affine connections (we write
$(M,p)$); the class is characterised by the fact that two connections
$\nabla$ and $\widehat\nabla$ in $p$ have the same geodesics up to
parametrisation. A model structure is $n$-dimensional projective space
$\mathbb{RP}^n$, but to avoid issues of orientability we prefer to
work with its double cover, the projective sphere $\mathbb{S}^n$. We
view this as a homogeneous space for $G=SL(n+1,\mathbb{R})$ and write
$P$ for isotropy subgroup of a point; so $P$ is a maximal parabolic
subgroup  and
we may identify $\mathbb{S}^n$ with $G/P$.  To each irreducible
$G$-representation $\mathbb{V}$ there is on $\mathbb{S}^n$ a canonical
finite resolution by linear differential operators \cite{CSS} which is
related to the algebraic resolutions from
\cite{BGG,Lepowsky},
\begin{equation}\label{BGG}
  0\to \mathbb{V}\to \cB^0 \stackrel{D}{\to} \cB^1\to\cdots \cB^n 
\to  0.
\end{equation} 
Here the $\cB^i$ are irreducible (weighted) tensor
bundles.  For the differential operators in the sequence there are
canonical {\em curved analogues}, that is, generalisations that exist
and are invariant on general projective manifolds $(M,p)$
\cite{CSSIII}. In particular, this is true for each {\em first BGG}
operator $D$, and we use the same terminology (and notation) for each
corresponding curved analogue.

These projectively invariant first BGG operators give the
equations we study. For these we consider the special class of
so-called normal solutions; see Section \ref{normS}. On projective
manifolds there is a canonical normal Cartan connection on a higher
frame bundle. The equivalent induced linear connections are termed (normal)
tractor connections \cite{BEG,CapGoTAMS}; the tractor and Cartan
connections are reviewed briefly in Section \ref{pT}. Beginning with a
tractor bundle induced from an irreducible $G$-representation, each
parallel tractor is equivalent to a normal solution of a first BGG
equation. On $\mathbb{S}^n$ all solutions are normal, but a
priori on curved structures it is a restriction.

By their definition, normal solutions are related to holonomy
reductions of the Cartan/tractor connection and from this perspective
certain local aspects have been investigated in
\cite{ArmstrongP1,ArmstrongP2}. There one sees that, on the one hand,
the available holonomy groups restrict the range of curved cases
(although the treatment there classifies irreducible holonomy
algebras, and so is not exhaustive), but on the other hand within the
allowed groups interesting geometric structures arise including
various pseudo-Riemannian Einstein (mentioned above) and contact
adapted projective structures in the sense of \cite{FoxIndiana}.

The main focus of the current article is to show how to access
geometric and topological information via a remarkable connection to
the model $\mathbb{S}^n$, and then a new understanding of the nature
of the first BGG solutions there.  Recalling that M\"{o}bius'
homogeneous coordinates are a fundamental tool for the calculus on
projective space, in Section \ref{chom} we find local curved
analogues, see Lemma \ref{ncoords3}. Based around a point, these are
determined uniquely by a normal tractor frame that we also build. The
whole construction is canonical up to a freedom parametrised by the
parabolic subgroup $P$.  The generalised homogeneous coordinates lead
to a diffeomorphism between the curved projective manifold $(M^n,p)$
and the model $\mathbb{S}^n$ that, together with the normal frame,
encodes a high degree of geometric contact.  In particular, in the
normal trivialisations, the components of parallel tractor fields on
$(M,p)$ pull back to parallel tractor fields on the model, see Theorem
\ref{canp1} and Corollary \ref{ccanp1} which are the first main
results.

An immediate consequence of these constructions is Corollary \ref{norms} which
shows that the components of normal solutions are necessarily the push
forward (via the diffeomorphism) of a first BGG solution on the model.
It follows that many local analytic and geometric questions for normal
solutions on $(M,p)$ can be settled by studying the same problem in
the simpler setting of the model. This strongly suggests there is
significant value in understanding the nature of first BGG solutions
on $\mathbb{S}^n$. It turns out that the answer, given in Proposition
\ref{BGG=Girr}, is rather appealing: the first BGG solutions on
$\mathbb{S}^n$ are precisely the (weighted) irreducible tensor fields
arising from a natural class (that we term $G$-irreducible) of
homogeneous polynomial systems on $\mathbb{R}^{n+1}$. For example, for
a given $k$, the solutions of \nn{scales}, on $\mathbb{S}^n$, are
simply the projective polynomials of degree $k$. Remarkably in the
curved setting, normal first BGG solutions are given by the same formal
polynomial systems, now understood as polynomials in the
generalised homogeneous coordinates.
 In general, these constructions allow us to conclude
that a normal solution $\tau$ of a first BGG operator is, in a precise
way, a curved analogue of a $G$-irreducible polynomial tensor field on
$\mathbb{S}^n$. This is an interpretation of the next main result, Theorem
\ref{canp}.

To indicate the scope we point out that appropriate collections of the
solutions of the BGG operators on scales \nn{scales} are sufficient to
yield curved analogues of {\em any} projective polynomial system see
Section \ref{parT}. Note that since we work over $\mathbb{R}$, a
polynomial system generally contains strictly more information that
the algebraic set it determines. In particular various distinct
geometric structures arise as curved analogues of distinct polynomial
polynomial systems with the common feature of empty zero locus; these
can be important and interesting, and in the analogous conformal
setting this includes Fefferman spaces \cite{CGI,CGF}, and positive
Riemannian Einstein metrics
\cite{Go-al}.

As an immediate application, in Corollary \ref{zeros} we see that the
local nature of the zero locus of  a
normal solution $\tau$ may be completely deduced from the data of the
corresponding algebraic set on $\mathbb{S}^n$.  For example we can use
this to descibe classes of cases where any zero locus of $\tau$ is
necessarily a smooth embedded submanifold.  In fact information that
is both finer and has global content is available. For a projective
manifold $(M,p)$ equipped with a normal solution $\tau$ we obtain a
decomposition, or more accurately stratification, of $(M,p)$ which
reflects the corresponding Bruhat-type orbit decomposition of the
model $\mathbb{S}^n$; this is termed a $P$-type decomposition. This
perspective should be useful in developing curved analogues of the
vector valued Poisson transforms (cf.\ for example \cite{Orsted}).
Information of a more analytic nature can also be deduced from the
model, and indeed in two of the examples of Section \ref{prelim} we
use this idea to show that the open $P$-types are geodesically
complete.
 
Finally we should point out that we have selected here, for
development in some detail, just part of a very general picture. It is
essentially clear that a direct analogue of our constructions is
possible for conformal geometry; the importance of parallel tractors
is more established in the conformal setting \cite{A} and examples
include the Fefferman space \cite{CGI,CGF,L} (as well as PE and almost
Einstein geometries). In fact via a different approach related results
can be established for all Cartan geometries, and thus in particular
for all parabolic geometries
\cite{CGH}.

We thank Robin Graham for pointing out his construction with Fefferman
of structures equivalent to the Klein-Einstein manifolds, which arise
among the examples of Section \ref{examples}.

\section{The curved analogue of projective polynomial systems}

Here we shall construct and exploit curved analogues of certain
projective polynomial systems.  We require some background to describe
the construction.

\subsection{Projective differential geometry and tractor calculus}
As mentioned above, we shall write
$\mathbb{S}^n:=\mathbb{P}_+(\mathbb{R}^{n+1})$ to denote the ray
projectivisation of $\mathbb{R}^{n+1}$.  This has a natural class of
preferred paths that may be viewed as unparametrised geodesics; these
arise from the projectivisation of 2-dimensional linear subspaces in
$\mathbb{R}^{n+1}$.  This structure is preserved by a group
action. Evidently, $G:=SL(n+1,\mathbb{R})$ acts transitively on
$\mathbb{S}^n$ and maps geodesics to geodesics.  To be concrete in our
development, we fix some choice of non-zero $e_0\in \mathbb{R}^{n+1}$
and define $P$ to be the parabolic subgroup stabilising the ray
$\mathbb{R}_+\cdot e_0$.

The classical curved generalisation of $\mathbb{S}^n$ is termed a {\em
  projective structure} $(M^n,p)$, $n\geq 2$, as defined in the introduction. 
Alternatively phrased, as connections on $T^*M$, the elements in $p$
satisfy
\begin{equation} \label{ptrans}
 \widehat\nabla_a u_b = \nabla_a u_b - \Upsilon_a u_b -\Upsilon_b u_a
\end{equation}
where $\Upsilon$ is some smooth section of $ T^*M$.

\subsection{Projective tractor calculus}\label{pT}

If $M$ is oriented we write $\ce(1)$ for the $(-n-1)$st root of the
canonical bundle. Otherwise we write $\ce(1)$ for a choice of line
bundle with $(-2n-2)$nd power the square of the canonical bundle.
We note that any
connection $\nabla \in p$ determines a connection on $\ce(1)$ and its
real powers $\ce(w)$, $w\in \mathbb{R}$; we call $\ce(w)$ the bundle
of projective densities of weight $w$. 
As a point on notation: Given a bundle $\mathcal{B}$ we shall write
$\mathcal{B}(w)$ as a shorthand for $\mathcal{B}\otimes \ce(w)$.

Although, by the definition of a projective structure, there is no
preferred connection on $TM$, there is a canonical connection, known
as the tractor connection, on a related higher
rank bundle. In the case of the model $\mathbb{S}^n=G/P$ this is
remarkably simple. For any $P$-representation $\mathbb{W}$ one has the
induced homogeneous bundle $\mathcal{W}:=G\times_P \mathbb{W}$ where
this means $G\times \mathbb{W}$ modulo the equivalence relation
$(gr,v)\sim (g,r\cdot v)$, for $r\in P$. However in the special case
that $\mathbb{W}$ is the restriction to $P$ of a $G$-representation
then $\mathcal{W}$ is canonically trivialised $\phi: \mathcal{W}\to
(G/P)\times \mathbb{W}$ by $(g,v)\mapsto (gP,g\cdot v)$. Canonically
we have the trivial connection on $(G/P)\times \mathbb{W}$ and via
$\phi$ this pulls back to the tractor connection $\nabla^\cT$ on
$\cW$.

Since it occupies little space, and because in any case we need the
notation and concepts, we review briefly the construction of the
tractor connection in general; we follow \cite{BEG} and the
conventions there.  In an abstract index notation let us write $\ce_A$
for $J^1\ce(1)$, the first jet prolongation of
$\ce(1)$.  Canonically we have the jet exact sequence
\begin{equation}\label{euler}
0\to \ce_a(1)\stackrel{Z_A{}^a}{\to} \ce_A \stackrel{X^A}{\to}\ce(1)\to 0,
\end{equation}
where we have written $X^A\in \Gamma \ce(1)$ for the jet projection,
and $Z_A{}^a$ for the map inserting $\ce_a(1)$; these are both
canonical.  We write $\ce_A=\ce_a(1)\rpl \ce(1)$ to summarise the
composition structure in \nn{euler}.  As mentioned, any connection
$\nabla \in p$ determines a connection on $\ce(1)$, and this is
precisely a splitting of \nn{euler}.  Thus given such a choice we have
the direct sum decomposition $\ce_A
\stackrel{\nabla}{=} \ce_a(1)\oplus \ce(1) $ with respect to which we
define a connection by
\begin{equation}\label{pconn}
\nabla^{\mathcal{T}}_a (\mu_b ~ \mid ~ \si)
:= (\nabla_a \mu_b + P_{ab} \si ~ \mid ~ \nabla_a \si -\mu_a). 
\end{equation}
Here $P_{ab}$ is the projective Schouten tensor and, with
$R_{ab}{}^c{}_d$ denoting the curvature of $\nabla$, is related to the
Ricci tensor $R_{ab}:=R_{ca}{}^c{}_b$ by
$(n-1)P_{ab}=R_{ab}-\frac{2}{n+1}R_{[ab]}$; $[\cdots]$ indicates the
skew part over the enclosed indices.  It turns out that \nn{pconn} is
independent of the choice $\nabla \in p$, and so
$\nabla^{\mathcal{T}}$ is determined canonically by the projective
structure $p$. This is the {\em cotractor connection} of \cite{Thomas}
and is equivalent to the normal Cartan connection for the Cartan
structure of type $(G,P)$, see \cite{CapGoTAMS}. Thus we shall also
term $\ce_A$ the {\em cotractor bundle}, and we note the dual {\em
  tractor bundle} $\ce^A$ (or in index free notation $\mathcal{T}$)
has canonically the dual {\em tractor connection}: in terms of a
splitting dual to that above this is given by
\begin{equation}\label{tconn}
\nabla^\cT_a \left( \begin{array}{c} \nu^b\\
\rho 
\end{array}\right) = 
\left( \begin{array}{c} \nabla_a\nu^b + \rho \delta^b_a\\
\nabla_a \rho - P_{ab}\nu^b
\end{array}\right). 
\end{equation}

It will be useful to understand how the underlying connections in $p$
arise from the tractor connection. By dualising \nn{euler} it follows
that the tractor bundle has a canonical composition structure given by
the exact sequence
\begin{equation}\label{ceuler}
0\to \ce(-1) \stackrel{X^A}{\to} \ce^A \stackrel{Z_A{}^a}{\to} \ce^a(-1) \to 0.
\end{equation}
The isomorphism $\ce_A \stackrel{\nabla}{=} \ce_a(1)\oplus \ce(1) $
determined by $\nabla\in p$ also splits \nn{ceuler} and is evidently
equivalent to a choice of section $Y_B\in \ce_B(-1)$  satisfying $X^BY_B=1$. Such a splitting is
equivalent to a {\em Weyl structure} (cf.\ \cite{CSbook}).  This
determines a bundle monomorphism,
\begin{equation}\label{Z}
Z^A{}_a: \ce^a(-1) \to \ce^A .
\end{equation}
Using this a connection $\nabla^Y_a$ on $\ce^a(-1)$ is then recovered
from $\nabla^\cT$ by the composition (on sections of $\ce^a(-1)$) of
$\nabla^\cT$ with the map \nn{Z} followed by the canonical map
$\ce^A\to \ce^a(-1)$ of \nn{ceuler}; this is evident from
\nn{tconn}. This then determines a connection on $TM$ that we may denote $\nabla^Y$. 

Since $\ce_A$ is $J^1\ce(1)$ we have the canonical universal 1-jet
differential operator $\mathbb{D}_A :\ce(1)\to \ce_A$, and from \nn{euler}
$X^A\mathbb{D}_A $ is the identity on $\ce(1)$. Thus any nonvanishing section
$\si\in \Gamma \ce(1)$ determines a special Weyl structure, termed a
{\em scale}, by taking $Y_A:=\si^{-1}\mathbb{D}_A \si$. In fact by considering
powers and roots of $\ce(1)$ one sees that $D_A$ generalises to an
invariant operator $\mathbb{D}_A:\ce(w)\to \ce_{A}(w-1)$, $w\in \mathbb{R}$,
known to Thomas \cite{BEG,Thomas} and we shall term any non-vanishing
section $\tau\in \Gamma\ce(w)$, $w\neq 0$, a scale (since we may take
$Y_A=\frac{1}{w}\tau^{-1}\mathbb{D}_A \tau$). We  write $\nabla^\tau$
for the affine connection in $p$ determined by a choice of scale $\tau$.

\subsection{Normal solutions}\label{normS}
Any vector bundle which is a tensor product of tensor powers of the
tractor and cotractor bundles, or a tensor part thereof, is termed a
tractor bundle. The structures which arise are handled efficiently by
appeal to a principal bundle picture as follows.

Following \cite{CapGoTAMS} we consider the bundle $\cG$ of adapted
frames for $\mathcal{T}$ which respect the filtration structure shown
in \nn{ceuler}.  This is a principal bundle with structure group $P$.
Then tautologically $\mathcal{T}$ is the associated bundle
$\cG\times_P \mathbb{R}^{n+1}$.  It is also straightforward to
recover, from the tractor connection, the unique Cartan connection
$\omega$ on $\cG$ from which the tractor connection is induced. It
follows that, given any representation $\mathbb{W}$ of $G$, we obtain a
tractor bundle $\mathcal{W}=\cG\times_P \mathbb{W}$ equipped with a
(linear) tractor connection induced from
$\omega$.  When we talk about tractor fields being {\em parallel} we
mean that they are covariantly constant with respect to this
 connection.  If $\mathbb{W}$ is an irreducible
$G$-representation then we say that the tractor bundle $ \mathcal{W}$
is {\em $G$-irreducible}. In this case there is a natural bundle map
$\Pi: \mathcal{W}\to \cB^0$, where $\cB^0$ is an irreducible weighted
tensor bundle, induced by the $P$-epimorphism from $\mathbb{W}$ to its
$P$-irreducible quotient.

\begin{proposition}\label{normp} \cite{CSS}
  Let $\mathcal{V}$ be a $G$-irreducible tractor bundle on $(M,p)$ and
  suppose that $I$ is a parallel section of $\mathcal{V}$. Then the
  bundle map $\Pi: \mathcal{V} \to \cB^0$ takes $I$ to a solution
  $\tau:=\Pi(I)$ of a first BGG operator
$$
D:\cB^0\to \cB^1.
$$
\end{proposition}

\noindent{\bf Definition:} We shall say that $\tau$, arising as in the
Proposition, is a {\em normal solution} (of the operator $D$). In the
following text we may use the term ``normal solution'' to mean the
normal solution for some first BGG operator $D$, without specifying
$D$.

\noindent{\bf Remark:} It is worth noting that in the case of the
model $\mathbb{S}^n$ all solutions arise this way. In fact in the
resolution \nn{BGG} the $G$-representation $\mathbb{V}$ may be
identified with the space of parallel tractors in the tractor bundle
associated to $\mathbb{V}$.

\subsection{The Thomas cone space $\aM$}\label{ambient} In view of the
canonical fibration $\ul{\pi}:\mathbb{R}^{n+1}\setminus 0\to
\mathbb{S}^n$ we may regard
$\mathbb{R}^{n+1}\setminus 0$ as a cone space over $
\mathbb{S}^n$. Here we recover the curved analogue of this (which was
probably known to T.Y.\ Thomas), see \cite{CSbook,FoxIndiana}.

In Section \ref{normS} above we mentioned that the Cartan connection
induces a canonical tractor connection on any associated bundle
$\cG\times_P \mathbb{W}$, where $\mathbb{W}$ is the restriction to $P$
of a $G$-representation. It is immediate from the equivariance
properties of $\omega$ that, more generally for any closed subgroup
$P_0\subset P$, we obtain a canonical connection on $\cG\times_{P_0}
\mathbb{W}$, a vector bundle over the fibrewise quotient $\cG/P_{0}$.
In particular, let us henceforth write $P_0$ to denote the subgroup of
$G$ fixing $e_0$ and define $\aM$ to be the quotient $\cG/P_0$, that
is, it is the total space
$$
\aM= \cG\times_P \mathbb{E}_+
$$
where $\mathbb{E}_+ $ is the $\mathbb{R}_+$-ray generated by $e_0$ in
$\mathbb{R}^{n+1}$; from \nn{ceuler} we see that it is equivalently
the total space of the ray-bundle $\ce(-1)_+$ (i.e.\ the subbundle of
positive rays in $\ce(-1)$). We write $\pi:\aM\to M$ for the canonical
bundle projection.

Now observe that, as $P_0$ representations, we have ${\frak
  g}/\frak{p}_0 \cong \mathbb{R}^{n+1}$.  From this there follow two
points. First by the last isomorphism, and that $\mathbb{R}^{n+1} $
may be considered as the restriction to $P_0$ of a $G$-representation
space, it follows that $\omega$ canonically induces a vector bundle
connection on $\cG\times_{P_0} {\frak g}/\frak{p}_0$. Second, by the
standard theory of Cartan connections, we also have canonically the
identification $\cG\times_{P_0} {\frak g}/\frak{p}_0 \cong T\aM$. From
the formula for the tractor connection (equivalently the normalisation
conditions of the normal Cartan connection) it follows that this
connection is Ricci-flat.  In summary.
\begin{proposition} \label{ambprop} The projective structure $(M,p)$
  determines a canonical Ricci-flat affine connection $\and$ on the
  manifold $\aM$.
\end{proposition}

  The canonical section $X^A$ corresponds to a section $\zeta^A$ of
  $T\aM$ which generates the $\mathbb{R}_+$ action on the fibres, and
  it is straightforward to verify that
\begin{equation}\label{andX}
\and_B\zeta^A=\delta^A_B.
\end{equation}

Associated bundles on $M$ arise from $P$-representations $\mathbb{U}$
as $\cG\times_P \mathbb{U}$. Sections are functions $u: \cG\to
\mathbb{U}$ which are $P$-equivariant in the sense that $u(g\cdot
r)=r^{-1}\cdot u (g)$. Since $P$ equivariance trivially implies
equivariance for any subgroup, it follows that, by restriction, such
sections lift immediately to sections of the corresponding bundle
$\cG\times_{P_0} \mathbb{U}$ over $\aM$. In particular using the
formula for the tractor connection from \cite[Section 2.5]{CapGoTAMS}
one sees immediately that, in the case that $\mathbb{U}$ is a
$G$-representation, parallel tractor fields on $M$ correspond in an
obvious way with parallel tensor fields on $\aM$. It follows that
arbitrary smooth sections of $\mathcal{T}$ (or $\cT^*$) correspond to
sections of $T\aM$ (resp.\ $T^*\aM$) that are in the null space of
$\zeta^A\and_A$, and so general (unweighted) tractor fields correspond
in an obvious way to tensor fields that are parallel in the directions
of the fibres of $\pi$. Using \eqref{andX} and that $\and$ is torsion
free, this means that a section of $\cT$ corresponds to a section of
$T\aM$ which is homogeneous of degree $-1$, with respect to the
principal $\Bbb R_+$--action.

Now any section of $\pi$ determines a splitting of \nn{ceuler} and so
a connection from $p$.  Note that a scale $\si$ determines a unique
section of $\ce(-1)_+$ and thus a section of $\pi$. It is
straightforward to verify \cite{CSbook} that the affine connection
$\nabla^\si$ that arises is related to the Thomas space connection
$\and$ as follows.
\begin{lemma}\label{uconn}
  Let $u,v\in \Gamma TM$ and $\si$ a scale viewed, as a section of
  $\pi:\aM\to M$. Then 
$$
\nabla_u v= \pi_*(\and_{\si_*u} \si_* v).
$$
\end{lemma}

Note that \nn{andX} implies that $\zeta^A\and_A \zeta^B=\zeta^B $ so each
fibre of $\pi$ agrees with the trace of a vertical geodesic.
  It follows that other geodesics remain
transverse to the fibres for all time, and project to regular curves
on $M$. It is an easy consequence of \nn{uconn} that these are geodesics from 
the class on $(M,p)$.

\subsection{Generalised homogeneous coordinates}\label{chom}

Here we shall show that given a point $q\in M$, and a choice of
adapted frame for $\mathcal{T}(1)_q$, we obtain an otherwise canonical
diffeomorphism between $(\aM,\and)$ and affine $\mathbb{R}^{n+1}$;
this map is distinguished by its properies of geometric contact with
the model, as we shall see later in this section.

Recall we denote by $\zeta$ the fundamental vector field generating
the principal right $\Bbb R_+$--action on $\aM$. In the case of the
model $\ul{\pi}:\mathbb{R}^{n+1}\setminus 0\to \mathbb{S}^n$
 the fundamental field coincides with the
usual Euler vector field $E$, and the affine connection $\and$ agrees
with the usual affine parallel transport.
\begin{lemma} \label{ncoords3} Choose $\tilde{q}\in \aM$, and a unit
  volume frame $e_0,\dots,e_n$ for $T_{\tilde{q}}\aM$, with
  $e_0=\zeta$.  This determines a diffeomorphism $\Phi: \ul{\pi}^{-1}
  (U') \to \pi^{-1}(U)$ for some open neighbourhood $U$ of
  $q:=\pi(\tilde{q})$ and some open
  set $U'$ in $\mathbb{S}^n$. With the following properties:\\
  $\small{\bullet}$ $\Phi$ is $\mathbb{R}_+$-equivariant and so
  determines
  a diffeomorphism $\phi:U'\to U$;\\
  $\small{\bullet}$ $\Phi$ maps straight lines through
  $\Phi^{-1}(\tilde{q})$ to geodesics for $\and$ through $\tilde{q}$,
  and so $\phi$ maps great circles through
  $\phi^{-1}(q)$ to geodesic paths through $q$;\\
  $\small{\bullet}$ $\Phi^*\zeta$ is the Euler vector field on
  $\ul{\pi}^{-1}(U')\subset\Bbb R^{n+1}$.
\end{lemma}

\noindent{\bf Proof:}
We shall write $\exp$ for the affine exponential map of
$\and$ at the point $\tilde q$. Now let $W$ be an open neighborhood of
zero in $\Bbb R^n$ such that 
$$
(x^1,\dots,x^n)\mapsto \pi(\exp(x^1e_1+\dots+x^ne_n)) 
$$
defines a diffeomorphism from $W$ onto an open neightborhood $U$ of
$q:=\pi(\tilde q)$ in $M$. We may identify $W$ with the affine
hyperplane neighbourhood $\{(1,x):x\in W\}$ in $\mathbb{R}^{n+1}$, and
write $U'\subset \mathbb{S}^n$ for the open subset consisting of its image under $\ul{\pi}$.
Now define a map
$\Phi:\underline{\pi}^{-1}(U')\to \pi^{-1}(U)$ by
$$
(r,rx^1,\dots,rx^n)\mapsto \exp(x^1e_1+\dots+x^ne_n)\cdot r,
$$
where $r>0$, $(x^1,\dots,x^n)\in W$ and the dot indicates the
principal right action. Evidently, this is an $\Bbb R_+$--equivariant
diffeomorphism, so it induces a diffeomorphism $\phi:U'\to U$ and
$\Phi^*\zeta$ is the Euler vector field on
$\ul{\pi}^{-1}(U')\subset\Bbb R^{n+1}$. Also, $\phi$ maps great
circles through $\ul{\pi}((1,0))$ to geodesic paths through $q$, since
by construction it maps straight lines in the affine hyperplane
through $(1,0)$ to geodesics through $\tilde q$ in $\aM$. 
\quad $\Box$

\smallskip

\noindent{\bf Remark:} Note that the frame $\{e_0,\dots,e_n\}$ for
$T_{\tilde q}\aM$ determines an adapted frame for
$\mathcal{T}_q$. Varying $\tilde{q}\in \pi^{-1}(q)$, any adapted frame
can be obtained in this way. Hence at a given point $q\in M$, the
freedom of choice is parametrised by $P$.

\smallskip

\noindent The Lemma leads to the following observation:\\
\noindent{\bf Remark:} {\em Generalised homogeneous coordinates.}
Let us write $X^0,X^1,\cdots ,X^n:\pi^{-1}(U)\to\Bbb R$ for the
functions on $\pi^{-1}(U)\subset \aM$ which are the push forward via
$\Phi$ (i.e.\  pull back via $\Phi^{-1}$) of the standard coordinates
$X^0,X^1,\cdots ,X^n$ on $\Bbb R^{n+1}$ (restricted to
$(\ul{\pi})^{-1}(U')$). Since $\Phi$ is a diffeomorphism, the
$X^{\ul{A}}$, $\mbox{\scriptsize{$\ul{A}$}}=0,1,\cdots ,n$, are coordinates on
$\pi^{-1}(U)$. Also note that by the equivariancy of $\Phi$, these
functions are homogeneous of degree one for the principal $\Bbb
R_+$--action on $\aM$, so they are equivalent to 1--densities on $M$.
This collection of densities may be viewed as curved versions of
homogeneous coordinates.

\subsection{The fundamental theorem for parallel
  tractors} \label{parT} We show here that the diffeomorphism of Lemma
\ref{ncoords3} captures a high degree of contact between $(M,p)$ and
$\mathbb{S}^n$. This is observed by a compatibility between the
tractor parallel transport, on the two manifolds, that we shall
describe precisely.  First we construct a frame field for $T\aM$ on
$\pi^{-1}(U)$ that corresponds to an adapted frame field for $\cT$ on
$U$.

Here we continue the notation of Lemma \ref{ncoords3}.
Take the vectors $e_1,\dots,e_n$ at $\tilde q$, and transport them
parallely along the horizontal geodesics $t\mapsto \exp(tx)$ for $x$
in the span of $e_1,\dots,e_n$. Possibly shrinking $U$, these vectors
project onto a local frame $\{\xi_1,\dots,\xi_n\}$ for the tangent
bundle $TM$ over $U$. 

Next, we claim that putting $e_0=\zeta$ along these horizontal
geodescis, we obtain a unit volume frame $\{e_0,\dots,e_n\}$ along
$\exp(W)$. Let $c(t)$ be one of the horizontal geodesics through
$\tilde q$. Then $c'(t)$ is obtained by parallely transporting $c'(0)$
along the geodesic to $c(t)$. By assumption, $c'(0)$ lies in the span
of $e_1,\dots,e_n$, whence $c'(t)$ lies in the span of
$e_1(c(t)),\dots,e_n(c(t))$. But together with $\and_\xi\zeta=\xi$ and
$\and_{c'(t)}e_i=0$ along $c(t)$, this implies that
$\and_{c'(t)}(\zeta\wedge e_1\wedge \dots\wedge e_n)=0$ along
$c(t)$. Since $\and$ is volume preserving, and the frame has unit
volume in $\tilde q$, the claim follows.

Finally we extend our frame along the filbres of $\pi$ by requiring
homogeneity of degree $-1$ with respect to the principal $\Bbb
R_+$--action, that is we require $e_i(y\cdot r)=r^{-1}\rho^r_* \cdot
e_i(y)$, where $\rho$ denotes the $\mathbb{R}_+$-action. Then it is
clear by construction that $e_i$ defines a frame for the tangent
bundle $T\aM$ over $\pi^{-1}(U)$, and at the same time determines an
adapted frame for $\cT$ over $U$ via the correspondence of Section
\ref{ambient}. Notice that by construction and the equivariancy of
$\Phi$ we see that since $\Phi^*e_0$ equals the Euler vector field $E$
on $\{(1,x):x\in W\}$, we have $\Phi^*e_0=(X^0)^{-1}E$ on
$\ul{\pi}^{-1}(U')$.

Next we need to known what this construction yields on the model
$\mathbb{S}^n$. The natural choice is to take $\tilde q=(1,0)\in
\mathbb{R}^{n+1}$ and
$e_{\ul{A}}=\partial_{\ul{A}}=\frac{\partial}{\partial
  X^{\ul{A}}}(\tilde q)$ for $\ul{A}=1,\dots,n$. It is easily
concluded that the construction just gives the frame field
$\{\frac1{X^0}E,\partial_1,\dots,\partial_n\}$ on the half space
$X^0>0$.  Now a constant tensor $\aI'$ on $\mathbb{R}^{n+1}$ is
equivalent to a parallel tractor on $I'$ on $ \mathbb{S}^n$; this uses
Section \ref{ambient} for the model.  Putting these things togther we
come to the following key fact. Here to simplify the statement, density
bundles are trivialised by the scale $X^0$ (corresponding to working on the section $X^0=1$ of $\aM$).
\begin{theorem}\label{canp1} 
  Suppose that $I$ is a parallel section of a tractor bundle.
  Composing with $\phi$ the coordinate functions of $I$ with respect
  to the frame derived from $\{e_0,\dots,e_n\}$, one obtains the
  coordinate functions of a parallel tractor $I'$ on the homogenous
  model with respect to the tractor frame obtained in the same way
  from $\{\frac1{X^0}E,\partial_1,\dots,\partial_n\}$.
\end{theorem}
\noindent{\bf Proof:}  Denote by $\aI$ the parallel tensor on  
$(\aM,\and)$ equivalent to $I$. We consider an
expression for $\aI$ in the form $a_J e_J$, where the
elements $e_J$ are linear combinations of tensor products of the
$e_j$, which form a local frame for the given tensor bundle. Then
along any of the geodesics $c(t)$ through $\tilde{q}$, and lying  in $\exp(W)$,
we can consider
$$
0=\and_{c'(t)}\sum_Ja_Je_J=\sum_J(c'(t)\cdot a_J(c(t)))e_J+\sum_J
a_J\and_{c'(t)}e_J 
$$
To expand the last term, we only need to know that
$\and_{c'(t)}e_0=c'(t)$ while $\and_{c'(t)}e_i=0$ for $i>0$ along
$c(t)$. This shows that for the coefficients $a_J$ we obtain a first
order ODE on the function $t\mapsto a_J(c(t))=b_J(t)$ which has the
form $b'(t)=F(b(t))$ and one obtains the same system on the corresponding straight line through $(1,0)\in\mathbb{R}^{n+1}$ over the homogenous model
$\mathbb{S}^n$. In vertical directions, everything is fixed by
homogeneity, so we obtain the result.  \quad $\Box$

There is a useful variant of the above result.  To simplify the
discussion let us simply identify tractor fields in $M$ with the
corresponding  homogeneous tensor fields on $\aM$
and do the same on the model.

On the homogeneous model, we can directly compute the change from the
frame $\{\frac1{X^0}E,\partial_1,\dots,\partial_n\}$ to
$\{\partial_0,\dots,\partial_n\}$, and then make the same change on
$M$. Denoting, as above, the generalised homogeneous coordinates on $M$ by
$X^{\ul{A}}$, this implies that any parallel tractor has constant
coordinate functions with respect to the the frame
\begin{equation}\label{frame}
\{f_0,f_1,\dots,f_n\}, \quad \mbox{where} \quad
f_0=e_0-\frac{X^1}{X^0}e_1-\dots-\frac{X^n}{X^0}e_n 
\end{equation}
and $f_i=e_i$ for $i=1,\cdots ,n$. The corresponding coframe is 
given by $f^0=e^0$
and $f^{i}=e^{i}+\frac{X^{i}}{X^0}e^0$ for $i=1,\cdots
,n$. 

Since the change of frame \nn{frame} is rational in the coordinates
and these coordinates are the push forward by $\Phi$ of the standard
coordinates on $\mathbb{R}^{n+1}$, Theorem \ref{canp1} can be
 equivalently phrased in terms of the frame
$f_{\ul{A}}$.  More generally we may define, on a sufficiently small
neighbourhood $U$ of any point $q\in M$, a map $\overline{\Phi}$ from
tractor fields on $U'=\phi^{-1}(U)\subset \mathbb{S}^n$ to tractor
fields on $U$ as follows: Use the frame $f_{\ul{A}}$ and its dual to
trivialise the tensor bundles on $\aM$. Use the standard
$\mathbb{R}^{n+1}$ frame to do the same on the model. Then push
forward the component functions of tractors fields (as homogeneous
functions on $\mathbb{R}^{n+1}$) on the model via the diffeomorphism
$\Phi$ and interpret as components of a tractor field in the
trivialisation on $M$.  Then we have the following.
\begin{corollary}\label{ccanp1} Given $q\in M$, a choice of adapted frame
  $e_{\ul{A}}(q)$ for $\cT_q$ canonically determines, for some
  neighbourhood $U$ of $q$, a diffeomorphism $\phi:U'\subset
  \mathbb{S}^n\to U$, a trivialisation of tractor and density bundles, and a
  compatible map $\overline{\Phi}$ from tractor fields on $U'$ to
  tractor fields on
  $U$ with the properties:\\
  $\small{\bullet}$ for parallel tractors $I$ the component functions
  are constant; \\
  $\small{\bullet}$ any parallel tractor $I$ on $U$ is the image under
  $\overline{\Phi}$ of a parallel tractor
  $U'\subset \mathbb{S}^n$;\\
  $\small{\bullet}$ the components of the canonical tractor field
  $X^A$ (in the trivialisation) are exactly the generalised
  homogeneous coordinates $X^{\ul{A}}$ of Section
  \ref{chom}, and these are the image under $\overline{\Phi}$ of the standard coordinates on $\mathbb{R}^{n+1}$.
\end{corollary}
\noindent{\bf Proof:} The first part of the last fact follows from
\nn{frame} and that $X^0e_0=\zeta$.  The final observation is
immediate since on homogeneous functions on $\mathbb{R}^{n+1}$,
$\overline\Phi$ is just $\Phi_*$.  The other points were treated
above. \quad $\Box$

\smallskip

\noindent{\bf Remark:} By the same argument that led to the
first bullet point of the Corollary, we see that the frame field
$f_{\ul{A}}$ is parallel along those geodesics through $\tilde{q}\in
\aM$ which lie in $\exp(W)$. Then in the vertical directions we have
$\and_\zeta f_{\ul{A}}=0$, 
$\mbox{\scriptsize{$\ul{A}$}}=0,1,\cdots,n$. These properties with
$f_{\ul{A}}(\tilde{q}):=e_{\ul{A}}(q)$ obviously characterise this
frame, which we will call a {\em normal frame}. Thus the generalised
homogeneous coordinates are characterised as the component functions
(densities) of the canonical tractor $X^A$, with respect to this
normal frame. This frame also determines a {\em normal scale}
$\si=f^0_AX^A=e^0_AX^A$. This agrees with $\overline{\Phi}(X^0)$, so
trivialising density bundles on $M$ using $\si$, is compatible via
$\phi$ with trivializing density bundles on $\mathbb{S}^n$ using
$X^0$. This is implicit in the construction proving Theorem
\ref{canp1}.

In our discussion below, we shall use both the normal frame
$f_{\ul{A}}$ and the adapted frame $e_{\ul{A}}$. Note that the map
$\overline{\Phi}$ can also been obtained using the trivialisations
corresponding to the adapted frame.

\smallskip

Theorem \ref{canp1} and its equivalent Corollary \ref{ccanp1} allow us
to treat normal solutions since, by Proposition \ref{normp}, each 
such arises as $\Pi(I)$ for a parallel tractor $I$, where $\Pi$
denotes the projection $\Pi$ to the irreducible quotient bundle. This
is easily understood using the adapted tractor frame
$\{e_0,\cdots,e_n\}$ from above.  Referring to \nn{euler} and
\nn{ceuler}, but using the normal scale $\si=e^0_AX^A$ to trivialise
densities,
the projection $Z_A{}^a:\ce^A\to \ce^a$ to the irreducible quotient is
characterised by $e_0\mapsto 0$ and $e_i\mapsto\xi_i$ for
$i=1,\dots,n$. Dually for $\ce_A$ the projection is given by
contracting with $X^A$, and so is characterised by by $e^0\mapsto 1$
and $e^i\mapsto 0$ for $i=1,\dots,n$. These determine $\Pi$ on tensor
products.

Let us describe the trivialisation of tensor bundles on $M$ (or
$\mathbb{S}^n$) induced by the frame $\{\xi_1,\dots,\xi_n\}$ and the
normal density $\si$ (the frame on $\mathbb{S}^n$ given by the
projection of $\{\partial_1,\dots,\partial_n\}$ and the normal density
$X^0$) as a {\em normal trivialisation}.

\begin{corollary}\label{norms}
  Suppose that $\tau$ is a normal solution of a first BGG equation
  $D\tau =0$ on $(M^n,p)$.  For each point $q\in M$ there are open
  sets $U\ni q$ and $U'\subset \mathbb{S}^n$ and a diffeomophism
  $\phi:U' \to U$, which maps great circles through $\phi^{-1}(q)$ to
  geodesics through $q$, such that $\ul{\tau}\o\phi = \ul{\tau}'$
  where $\tau'$ is a solution for $D$ on the model, and $\ul{\tau}$, $
  \ul{\tau}'$ are the component functions of $\tau$ and $\tau'$ in the
  appropriate normal trivialisations. 
\end{corollary}

The Corollary here is a first version of what we shall shortly refer
to as the fundamental theorem for normal solutions. What is lacking at
this point is some explicit understanding of what the normal BGG
solutions are on the model.  Fortunately the tools we have already
developed give an answer almost immediately. Let us follow through a
the above for a special class of examples. 

Consider the case of a completely symmetric parallel tractor
$H_{A_1\cdots A_k}$ on $(M,p)$.  From \nn{euler} the map $\Pi$ of
Proposition \ref{normp} is simply $H_{A_1\cdots A_k}\mapsto
H_{A_1\cdots A_k}X^{A_1}\cdots X^{A_k}=:\si \in \Gamma\ce(k)$.  Now
$\si$ corresponds to degree $k$ homogeneous function on $\aM$ that we
shall denote the same way. From Corollary \ref{norms}, locally (and
using the constructions and notation from above) we have $\si = \phi_*
\si'$ where $\si'$ is a solution to $D$ on $\mathbb{S}^n$ .  With
respect to the trivialistions arising from the tractor frame
$f_{\ul{A}}$, we have $\si$ corresponds to the homogeneous function
$\si = H_{\ul{A}_1\cdots \ul{A}_k}X^{\ul{A}_1}\cdots X^{\ul{A}_k}$ on
$\aM$. From Corollary \ref{ccanp1}, the components $H_{\ul{A}_1\cdots
  \ul{A}_k}$ are constant, so $\si$ is expressed as a polynomial.
According to the construction, and using Corollary \ref{canp1}, $\si'$
is simply {\em the same} homogeneous polynomial (in standard
coordinates), now viewed as a homogeneous function on
$\mathbb{R}^{n+1}$. The latter is precisely a projective polynomial on
$\mathbb{S}^n$ of degree $k$.
 
Note that this outcome is also clear via another perspective from
Corollary \ref{ccanp1}: $\si=\Pi(H)$ is $\overline{\Phi}\si'$, for
some density $\si'$ on $\mathbb{S}^n$. By expressing $\si'$ in the
frame $f_{\ul{A}}$ as above, $\si= H_{\ul{A}_1\cdots
  \ul{A}_k}X^{\ul{A}_1}\cdots X^{\ul{A}_k}$, and using Corollary
\ref{ccanp1}, we see $\si'$ is given by the same formal expression on
$\mathbb{R}^{n+1}$ and hence is $\Pi(H')$ for a parallel tractor $H'$
there.  From either persepective, as a special case of
this we may replace $(M,p)$ with $\mathbb{S}^n$ to obtain the
following.
\begin{proposition}\label{symm}
  On $\mathbb{S}^n$ the degree $k$ projective polynomials are precisely the 
first BGG solutions corresponding to symmetric rank $k$ parallel
  tractors.
\end{proposition}
So we see that, at least for this class of cases, first BGG solutions
on the model are just projective polynomials. On the other hand
Corollary \ref{norms} shows that corresponding normal solutions in
$\Gamma\ce(k)$, on a projective manifold $(M,p)$, are curved
analogues.  Since any algebraic set arises from a collection of such
polynomials, we have a universal way to describe a canonical curved
analogue of the projective polynomial system involved.

We may elaborate on the Proposition somewhat. The PDE involved are
the equations of the operators $D_k:\Gamma\ce(k)\to
\Gamma\ce_{(a_1a_2\cdots a_{k})}(k)$ 
that, in terms of $\nabla\in
p$, take the form \nn{scales}.
The equation $D_k \si =0$ {\em characterises} degree $k$ projective
polynomial densities. The specific polynomial densities solving this
are then in 1-1 correspondence with parallel symmetric cotractors in
$\Gamma\ce_{(A_1\cdots A_k)}$. The parallel tractor corresponding to a
particular solution $\si$ is really part of the jet (roughly, Taylor
series) data of $\si$ determined by {\em prolongation} (in \cite{CSS}
and cf.\ \cite{BCEG}), but remarkably we can avoid any significant
details of this beyond what is implicit in the treament above.

We now see a generalisation of the Proposition above and corresponding
refinement of Corollary \ref{norms}. A first observation is that the
tractor connection preserves a section of $\Lambda^{n+1}\ce_A$ (a
tractor volume form), and it follows that, without loss of generality,
we may work with covariant tractor fields, and we henceforth make this
simplification.  However to work with these, and state results
concisely, we shall need to recall some standard representation
theory, as well as related notions and notation.

\subsection{Tensors and representations}

Up to isomorphism, each irreducible representation
of $G$ is described by a weight
\begin{equation}\label{young}
\underline{r}=(r_1,\cdots,r_n), \quad \mbox{where} \quad  r_1\geq \cdots \geq r_n\geq 0,
\end{equation}
and $r_i \in \mathbb{Z}_{\geq 0}, i=1,\cdots ,n$.  Below we use the
notation $|\underline{r}|:=\sum_1^n r_i$.  Equivalently the
representation \nn{young} is given by a Young diagram where from the
top, and proceeding down, the rows have respective lengths
$r_1,\cdots,r_n$, see e.g.\ \cite{FultonH}.  As a shorthand for
the weights we shall omit any terminal string of $0$s. For exmple
$(2,2)$ means
$$
\underbrace{(2,2,0, \cdots,0)}_{n}.
$$

\noindent{\bf Definition.} $G$-type and $P$-type. Let us view
$(\mathbb{R}^{n+1})^*$ as the dual of the defining representation of
$SL(n+1,\mathbb{R})$.  $G$ and its subgroups act on the tensor algebra
of $(\mathbb{R}^{n+1})^*$ and tensors will be said to be of the same
{\em $G$-type} (respectively $P$-type) if they lie in the same
$G$-orbit (respectively $P$-orbit).

\medskip

Each tensor power of $\bigotimes^{|\ul{r}|} (\mathbb{R}^{n+1})^*$ may
be decomposed into irreducible
representations classified by the weights as in \nn{young}.  Realising
the irreducible representations in tensor powers of the (dual to the)
standard representation is not unique. Certainly a tensor belonging to
the representation \nn{young} has valence (i.e.\ total rank)
$|\underline{r}|:=\sum_1^n r_i$.  The representation given in
\nn{young} will be realised by tensors
$$
I_{A_1\cdots A_{r_1}B_1\cdots B_{r_2}\cdots E_1\cdots E_{r_n}}
$$
which are completely symmetric over each of the respective index
sets $A_1\cdots A_{r_1}$, $B_1\cdots B_{r_2}$, and so on to $E_1\cdots
E_{r_n}$.  We call a tensor in such a subrepresentation
$G$-irreducible.  The irreducibility of the representation is further
encoded in what are sometimes termed hidden symmetries of the tensor
elements \cite{ot}.  For example symmetrising over any $r_1+1$ indices
will annihilate the tensor $I$.  With this understood we shall write
$\mathbb{R}(\ul{r})$ for the vector space of such tensors in
$\bigotimes^{|\ul{r}|} (\mathbb{R}^{n+1})^*$.

Finally as a point of notation. Above we have expressed the tensor $I$
as an object adorned with abstract indices. $(\mathbb{R}^{n+1})^*$ has
a standard basis, this generates a standard vectorial basis for the
tensor algebra it generates. In terms of this we may express $I$ in
terms of its components, and write
$$
I_{\ul{A}_1\cdots \ul{A}_{r_1}\ul{B}_1\cdots \ul{B}_{r_2}\cdots
  \ul{E}_1\cdots \ul{E}_{r_n}}.
$$

\subsection{$G$-irreducible polynomial systems}

As above, let us write $X^{\ul{A}}$, 
$\ul{\mbox{\scriptsize$A$}}=0,1,\cdots,n$ for the standard
coordinates on $\mathbb{R}^{n+1}$. Given a tensor in
$\bigotimes^{|\ul{r}|} (\mathbb{R}^{n+1})^*$ we may construct
polynomial systems by contraction in the obvious way. For example if
$R_{A_1A_2B_1B_2}\in \mathbb{R}(2,2)$ then we may form the two (in general
non-trivial) polynomial systems
$$
P_{\ul{A}_2\ul{B}_1\ul{B}_2}=R_{\ul{A}_1\ul{A}_2\ul{B}_1\ul{B}_2}X^{\ul{A}_1}
\quad \mbox{and} \quad 
Q_{\ul{B}_1\ul{B}_2}=R_{\ul{A}_1\ul{A}_2\ul{B}_1\ul{B}_2}X^{\ul{A}_1}X^{\ul{A}_2}
$$
where repeated indices are summed (according to the Einstein summation
convention).  We shall term the latter of these {\em saturated} since
any contraction of $X^{\ul{C}}$ into $Q_{\ul{A}\ul{B}}$ will result in
anihilation, as a result of the (hidden) symmetries of $R$.  Such
polynomials form a natural class for many purposes, as shall be clear
shortly.

\medskip

\noindent{\bf Definition:} We shall say that a polynomial system is {\em
  $G$-irreducible} if it arises as
\begin{equation}\label{Girr}
Q_{\ul{B}_1\cdots \ul{B}_{r_2}\cdots
  \ul{E}_1\cdots \ul{E}_{r_n}}
=I_{\ul{A}_1\cdots \ul{A}_{r_1}\ul{B}_1\cdots \ul{B}_{r_2}\cdots
  \ul{E}_1\cdots \ul{E}_{r_n}} X^{\ul{A}_1} \cdots X^{\ul{A}_{r_1}}
\end{equation}
for some tensor $I\in \mathbb{R}(\ul{r})$. Note that the system is
homogeneous of degree $r_1$, and defines a projective algebraic set on
$\mathbb{S}^n$.

The system $Q$ given in \nn{Girr} is saturated and homogeneous. By
construction, as a tensor field on $\mathbb{R}^{n+1}$ it has
symmetries consistent with the representation $\ul{r'}=(r_2,\cdots,
r_n)$.  It now follows that these collectively imply that it
corresponds to a certain field $\ul{\tau}$ on $\mathbb{S}^n$. In fact
from \nn{euler} and the relationship between tractors bundle sections
and cone tensor fields, as described in Section \ref{ambient}, we
obtain that this tensor $\ul{\tau}$ is a section of weighted
irreducible tensor bundle $\ce(r_2,\cdots,r_n)(k)$, where $k=
|\ul{r}|$. Here $ \ce(r_2,\cdots,r_n)$ is the bundle of covariant
tensors having the Young symmetry $(r_2,\cdots ,r_n)$ and
$\ce(r_2,\cdots,r_n)(k)= \ce(r_2,\cdots,r_n)\otimes \ce(k)$.  We shall
say that a tensor field on $\mathbb{S}^n$ that arises in this way is a
{\em $G$-irreducible polynomial tensor field}.  The critical point
here is that we can give a precise differential characterisation of
these, that we shall come to now.

By their definition (and with Proposition \ref{normp}) we see the
weighted tensor bundles $\ce(r_2,\cdots,r_n)(k)$, $k \geq r_2+|r'| $
in $\mathbb{Z}$, are exactly the bundles $\mathcal B^0$ in the BGG 
complexes mentioned in the introduction (see \nn{BGG}).  Now we can
state the full extension of Proposition \ref{symm}.
\begin{proposition}\label{BGG=Girr}
  The $G$-irreducible polynomial tensor fields on $\mathbb{S}^n$ are
  precisely the solutions of first BGG operators $D:\cB^0\to \cB^1$:
  If $\cB^0$ is realised as the irreducible weighted tensor bundle
  $\ce(r_2,\cdots,r_n)(k)$ and a section $\tau$ thereof is a solution
  (i.e.\ $D\tau=0$), then the homogeneous polynomial system
  corresponding to $\tau$ is the saturate of some $I\in
  \mathbb{R}(r_1,\cdots,r_n)$ where $r_1=k-|\ul{r}'|$.
\end{proposition}
\noindent This result is an easy consequence of the main Theorem \ref{canp},
which follows shortly.

\noindent{\bf Remark:} Note that on $\mathbb{S}^n$ the homogeneous
coordinates $X^0,\cdots ,X^n$ linearly generate the full solution
space of the $k=1$ system \nn{scales}. Thus via the Proposition, these
polynomially generate all first BGG solutions. The Proposition also
gives the specific polynomial systems involved.

\subsection{The fundamental theorem of normal solutions}

\medskip

As observed above each domain bundle $\cB^0$, for a first BGG operator,
may be realised in the form $\ce(r_2,\cdots,r_n)(k)$. 
Let $\tau\in \Gamma \ce(r_2,\cdots,r_n)(k)$ be a normal solution.
Then $\tau=\Pi(I^\tau)$ where
$I^\tau$ is a $G$-irreducible parallel tractor field.
It is an easy consequence of the
filtration structure arising from \nn{euler} that $\Pi$ may be
realised explicitly by saturating $I^\tau$ with $X$ to yield
$$
Q^\tau_{B_1\cdots B_{r_2}\cdots
  E_1\cdots E_{r_n}}=I^\tau_{A_1\cdots A_{r_1}B_1\cdots B_{r_2}\cdots
  E_1\cdots E_{r_n}} X^{A_1} \cdots X^{A_{r_1}},
$$ 
(where $r_1=k-|\ul{r'}|$) and by then contracting the projectors
$Z^A{}_a$ (of \nn{Z}) onto all remaining indices. Note that $Q^\tau$
takes values in an irreducible subbundle of the weighted tractor
bundle $\ce_{B_1\cdots B_{r_2}\cdots E_1\cdots E_{r_n}}(r_1)$ and the
final step of contracting with the concatenation of $Z^A{}_a$
projectors is simply realising the isomorphism between this subbundle
and the weighted tensor bundle $\ce(r_2,\cdots,r_n)(k)$. Thus we
henceforth identify $\ce(r_2,\cdots,r_n)(k)$ with this subbundle (as
done implicitly in Corollary \ref{norms}) and thus 
$\tau$ and $Q^\tau$ are also to be identified.

Now fix an arbitrary point $\tilde{q}\in \aM$. In a neighbourhood of
$\tilde{q}$, and in terms of a normal frame field (as defined in
\nn{frame}, and see the Remark below that) we have
$$
Q^\tau_{\ul{B}_1\cdots \ul{B}_{r_2}\cdots \ul{E}_1\cdots
  \ul{E}_{r_n}} =I^\tau_{\ul{A}_1\cdots \ul{A}_{r_1}\ul{B}_1\cdots
  \ul{B}_{r_2}\cdots \ul{E}_1\cdots \ul{E}_{r_n}} X^{\ul{A}_1} \cdots
X^{\ul{A}_{r_1}}
$$
and, by Corollary \ref{ccanp1}, the components $
I^\tau_{\ul{A}_1\cdots \ul{A}_{r_0}\ul{B}_1\cdots \ul{B}_{r_1}\cdots
  \ul{E}_1\cdots \ul{E}_{r_n}}$ are constant. Using again Corollary
\ref{ccanp1} to interpret this on $\mathbb{S}^n$ we have the following.
\begin{theorem}\label{canp}
  Suppose that $\tau \in \Gamma \ce(r_2,\cdots,r_n)(k)$ is a normal
  solution of the equation $D\tau =0$ on $(M^n,p)$. For an arbitrary
  point $q\in M$, fix a an adapted frame at $q$. This determines a
  normal frame $f_{\ul{A}}$ (as in \ref{frame}), a local
  diffeomorphism $\phi:\mathbb{S}^n\to M$, and corresponding
  generalised homogeneous coordinates $X^{\ul{A}}$ in a neighbourhood
  of $q$.
  In terms of these the following hold:\\
  $\small{\bullet}$ With respect to the trivialisation of
  $\ce_{B_1\cdots B_{r_2}\cdots E_1\cdots E_{r_n}}(r_1)$, determined by
  $f_{\ul{A}}$, the coordinate functions $\ul{\tau}$ (of $\tau$) form
  a homogeneous polynomial system in the generalised
  homogeneous coordinates. \\
  $\small{\bullet}$ The collection $\Phi^* \ul{\tau}$ is given by the
  same formal polynomial system, where the $X^{\ul{A}}$ are now
  interpreted as the standard coordinates of $\mathbb{R}^{n+1}$.  With
  respect to the standard frame on $\mathbb{R}^{n+1}$, the collection
  $\Phi^* \ul{\tau}$ are the coordinate functions of a solution
  $\tau'$ of the equation $D\tau' =0$ on $\mathbb{S}^n$.
\end{theorem}

\noindent{\bf Remark:} 
Note that the polynomial system $\ul{\tau}$, as in the Theorem,
satisfies polynomial relations. These arise in an obvious way from the
fact that tractor section, equivalent to $\tau$, is saturated with
respect to contraction with $X^A$. For example in the case $\tau \in
\Gamma \ce_a(2)$
then the system $\ul{\tau}$ consists of the $n+1$ linear polynomials
$K_{\ul{A}\ul{B}}X^{\ul{B}}$ where the matrix of components of $K$ is
skew, i.e.\ $K_{\ul{A}\ul{B}}=-K_{\ul{B}\ul{A}}$. Thus there is the
one polynomial relation $K_{\ul{A}\ul{B}}X^{\ul{A}}X^{\ul{B}}=0$.

\smallskip

Next note that since the reasoning in the first part of the proof above
applies, in particular, when we begin with $\tau$ a solution of $D$ in
the ``flat case'' (i.e.\ on $\mathbb{S}^n$ with its standard
projective structure) as a corollary we have at once the Proposition
\ref{BGG=Girr}.

\subsection{The zero locus of normal solution } \label{zeroS}
The Theorem above is local in nature but it (or equivalently Corollary
\ref{norms}) has a global interpretation. Before we come to this we
need some simple observations, and a definition.

\noindent{\bf Definition.} {\em $G$-type and $P$-type of a point.}  Consider
a projective manifold $(M,p)$ equipped a with a normal solution
$\tau$, and let $I^\tau$ be the parallel tractor such that
$\tau=\Pi(I^\tau)$. For a point $q\in M$ we can choose a tractor frame
$e_{\ul{A}}$ for the tractor space $\cT_q$ at $q$, which is adapted in
the sense that $e_0$ is parallel to $X^A$. Using this frame, the
components of $I^\tau(q)$ define an element in the tensor algebra of
$(\mathbb{R}^{n+1})^*$. As we have noted in \ref{normS} any two such
frames are related by the action of an element of $P$. Hence the
$P$--orbit of this element, which we call the {\em $P$-type of the
  point $q$}, depends only on $(M,p)$, $\tau$, and $q$, and not on
further choices.  Of course, this implies that also its $G$--orbit is
well defined, and we call this the {\em $G$-type of the point $q$}

\smallskip 

\noindent Since a parallel tractor has constant components with
respect to the normal frames $f_{\ul{A}}$ which,
although not adapted, are volume mormalised, we obtain the following. 
\begin{theorem}\label{gtype}
  If $(M,p)$ is connected then any parallel tractor field $I$ 
  has a constant $G$-type.
\end{theorem}
\noindent An analogous statement is not available for $P$-type,
because an adapted frame is not parallel along any curve; this is clear
from \nn{andX}.

It follows from Theorem \ref{gtype} that, for a normal solution $\tau$ on
a connected manifold $(M,p)$, we may associate a \ul{fixed}
$G$-irreducible polynomial tensor field $\tau'$ on the model
$\mathbb{S}^n$, and this is obviously unique (up to a $G$-action that we
shall ignore). In this case we shall say that
$(\mathbb{S}^n,\tau')$ is {\em the model} for $(M,p,\tau)$
(meaning $(M,p)$ equipped with a normal solution $\tau$).

\smallskip

For a given normal solution $\tau$ on a
projective manifold $(M,p)$, its zero locus $\mathcal{Z}(\tau) $ is
not necessarily smooth.  However, given $q\in \mathcal{Z}(\tau)$,
Corollary \ref{norms} shows that there is a point $q'\in
\mathcal{Z}(\tau')\subset\mathbb{S}^n$ and a local diffeomorphism
from an open neighborhood of $q$ to an open neighborhood of $q'$ in
the model, which is compatible with resepect to the zero sets. In
particular, $\mathcal{Z}(\tau)$ cannot have worse singularities than
$\mathcal{Z}(\tau')$.
For example if $q'$ is a smooth point then, in a neighbourhood of $q$,
$ \mathcal{Z}(\tau)$ is an embedded smooth submanifold.  Hence the
problem of the classification of zero locus singularities for normal
solutions is reduced to a problem in real algebraic geometry. 
Note that it can be that the model algebraic set
is not (globally) smooth and yet the zero locus $\mathcal{Z}(\tau)$ is
smooth and embedded; indeed as an extreme case $\mathcal{Z}(\tau)$ may
be empty.

\smallskip

For emphasis we collect some of these definitions and consequences of
Theorem \ref{canp} into a statement.
\begin{corollary}\label{zeros} Let $\tau$ be a normal solution on a
  connected manifold $(M,p)$ and let $(\mathbb{S}^n,\tau')$ be the
  corresponding model.  If $q\in \mathcal{Z}(\tau)$ then there some
  $q'\in \mathcal{Z}(\tau')$ with the same singularity type.  In
  particular, if $\mathcal{Z}(\tau')$ is a smooth algebraic set
  then $\mathcal{Z}(\tau)$ is a smooth embedded submanifold.
\end{corollary}

\subsection{Orbit type decomposition of $M$} Here we observe that a
normal solution $\tau$ (equivalently a parallel tractor field $I^\tau$) 
determines a canonical stratification of $M$. 

Different points on the manifold may have the same $P$-type and this
establishes an equivalence relation for the points of $M$. Thus the
points of the manifold $M$ are partitioned according to $P$-type. On a
given structure $(M,p,\tau)$ there can be many orbit types.  For
example a point where $\tau$ vanishes is certainly in a different
$P$-type to a point where $\tau$ is not zero. In general the $P$-type
decomposition exposes considerably finer structure than this example
illustrates. We treat some examples in Section \ref{prelim} below, but
a more detailed analysis is the subject of \cite{CGH}.

Now we claim that the diffeomorphism $\phi$, of Lemma \ref{ncoords3},
preserves $P$-type. Precisely, the $P$--type of $I^\tau$ at $\phi(y)$
coincides with the $P$--type of $I^{\tau'}$ at $y$. To see this,
arguing as in the proof of Theorem \ref{canp1}, simply observe that
these tractors have the same coordinate expression in different
frames, both of which are adapted. This exactly means same they have
the same $P$--type.

Thus we may take the alternative view that the
maps $\phi$ transfers the $P$-type decomposition of the model
$(\mathbb{S}^n, \tau')$ onto the corresponding curved structure
$(M,p,\tau)$.  Part of the power of this lies in the following result
(discussed in greater detail in \cite{CGH}).
\begin{proposition}\label{PGI} Let $G_{I'}\subset G$ be the isotropy
  subgroup of a tensor $I'$ in $\bigotimes \mathbb{R}^*$. Viewing $I'$
  as a parallel tractor on $\mathbb{S}^n$, the $P$-type decomposition
  of $(\mathbb{S}^n,I')$ is the same as the orbit decomposition of
  $\mathbb{S}^n$ under the action of $G_{I'}$.
\end{proposition}
\noindent{\bf Proof:}
It is easily verified directly, that in either case the orbits are
naturally parametrised by the points of the double coset space $G_{I'}
\backslash G \slash P$.  \quad $\Box$

\smallskip

In general a projective manifold $ (M,p)$, with a parallel tractor
$I$, admits no action by $G_{I'}$, and so a $G_{I'}$-orbit
decomposition of $M$ makes no sense. However we have the simple but
surprising outcome that the $P$-type decomposition is well defined,
and this echoes the $G_{I'}$-orbit decomposition of the model.  For
example it follows from the Proposition \ref{PGI} that the $P$-types
of $ \mathbb{S}^n$ are non-self-intersecting smoothly immersed (in
fact initial \cite[Theorem 5.14]{KMS}) submanifolds and hence, via the
diffeomorphism $\phi$ of Lemma \ref{ncoords3}, we may conclude that so
are the $P$-types of $(M,p,I)$. Thus the $P$-types give a smooth
stratification of $M$. In summary we have the following.
\begin{theorem}\label{Ptype}
  Let $(M,p)$ be projective manifold equipped with a parallel tractor
  $I$.  Then $(M,p)$ is stratified according to a $P$-type
  decomposition.  The diffeomorphism $\phi$ of Lemma \ref{ncoords3}
  preserves the $P$-type and, in particular, the occurring $P$-types
  are a selection of those arising in the model $(\mathbb{S}^n,I')$.
\end{theorem}

\subsection{Geometry}

The presence of a parallel tractor equips $(M,p)$ with additional
geometric structure. Of course a parallel tractor may be reinterpreted
as a reduction of the (projective) tractor holonomy. By definition
this is additional geometric structure, and as mentioned in the
introduction, local and generic aspects of this have been explored in,
for example, \cite{ArmstrongP1,ArmstrongP2}. In the examples of the
next section we shall see, from our current point of view, how such
classical structures arise. In one case we shall see the structure is
necessarily non-Ricci-flat pseudo-Riemannian Einstein on open subsets.

More importantly for the directions here, we show that the parallel
tractor, along with Theorems \ref{canp1} and \ref{canp}, provides a tool
which can relate the geometry of a normal solution zero locus to that
on the complementary space.

\section{Examples}\label{examples}

\subsection{Preliminary observations} \label{prelim} 
In several of the examples below we consider normal solutions $\tau$
which are sections of a density bundle $\ce(w)$, where $w\neq 0$, and
the open $P$-types are submanifolds on which $\tau$ is
non-vanishing. On such a $P$-type, denoted by $M_+$ say, $\tau$ is a
scale and so naturally determines an affine connection $\nabla^\tau$
from the projective class (restricted to $M_+$) as discussed in
Section \ref{pT}. We want to explain how our method can be used to
prove results related to geodesic completeness of $\nabla^\tau$.
Consider a geodesic path in $M$ which leaves $M_+$, i.e.~which
intersects both $M_+$ and $\mathcal{Z}(\tau)$. Take a point
$q\in\mathcal{Z}(\tau)$ which lies in the closure of $M_+$ and carry
out the construction from Lemma \ref{ncoords3} for some point $\tilde
q$ over $q$. Then our geodesic path will become a geodesic $\gamma$
for the connection determined by the normal scale $\sigma$ determined
by the construction. 

Now take a point $x\in M_+$ which lies on $\gamma$ (and in the range
of the diffeomorphism $\phi$). Starting from $x$ and moving along
$\gamma$ in the direction of $q$, the point $q$ will of course be
reached in finite time. Now some reparametrization
$\hat\gamma$ of $\gamma$ will be a geodesic for $\nabla^\tau$, and it
may happen that this reparametrization has the effect that $q$ is no longer
reached in finite time.

If we assume that the original manifold $M$ is closed
then the only way for $M_+$ to be geodescially incomplete is that
geodesics leave $M_+$ in finite time. Hence if in the above
considerations $q$ is never reached in finite time, geodesic
completeness follows. 

The reparametrization from $\gamma$ to $\hat\gamma$ can be obtained as
the solution of a an ODE which depends only on the function describing
$\tau$ in the trivialization of the density bundle determined by
$\sigma$. But now it follows from the construction that the
diffeomorphism $\phi$ relates $\sigma$ and the function describing
$\tau$ in the trivialization determined by $\sigma$ to their
counterparts on the homogeneous model. Consequently, the
reparametrizations on the homogeneous model and on the curved manifold
are determined by the same ODE and hence coincide. Thus, if the point
$\phi^{-1}(q)$ on the homogeneous model is not reached in finite time
after the reparametrisation, then the same is true for $q$. In particular,
geodesic completeness carries over from the homogenous model to curved
geometries on closed manifolds.  

If in such a situation $\mathcal{Z}(\tau)$ is the boundary of $M_+$
(as in first two examples below) then it reasonable to call
$\mathcal{Z}(\tau)$ the {\em projective infinity} for
$(M_+,\nabla^\tau)$ (as an analogue of the term conformal infinity).
This is not meant to imply that that $\mathcal{Z}(\tau)$ necessarily
has a canonical projective structure.

\subsection{The parallel standard cotractor -- projective almost
  Ricci-flat} \label{RF} 

Here we obtain a structure which generalises that of an affine Ricci
flat manifold (cf.\ \cite{Rthesis}); the result is affine Ricci-flat
on an open dense set.  From the point of view of compactifications, it
yields a structure that is a compactification of a Ricci-flat manifold
that is a curved generalisation of the usual projective
compactification of the affine plane to a hemisphere via central
projection.  Thus, in particular, it is different to a conformal
compactification of such a space (which for the case of the Euclidean
plane is a 1-point compactification).

The structure in this case in a projective manifold 
 $(M,p)$ equipped with a parallel section $I_A$ of the
standard cotractor bundle. Let $\si$ denote $I_A X^A$. To find the
first BGG equation that this satisfies we may calculate with respect
to $\nabla\in p $.  We have $I_B\stackrel{\nabla}{=} (\mu_b~~\si)$ and
using the formula for the tractor connection \nn{pconn}, we see that
$I_B$ parallel implies that $\mu_b=\nabla_b\si$ and then
\begin{equation}\label{basic}
\nabla_a\nabla_b \si +P_{ab}\si=0; 
\end{equation}
this is the $k=1$ case of \nn{scales}.  All solutions of \nn{basic}
arise this way (i.e. any solution of \nn{basic} is normal), indeed a
prolongation of this equation determines the tractor connection
\cite{BEG}.

Where $\si$ is non-vanishing \nn{basic} is the equation that $\nabla$
is projectively Ricci-flat; more precisely in an open neighbourhood
with $\si$ nowhere zero the connection $\widehat{\nabla}$,
characterised by $\widehat{\nabla}\si=0$ and hence related to $\nabla$ via \nn{ptrans} with $\Upsilon_a$ given
by $\si^{-1} \nabla_a \si$, is Ricci flat. Thus the structure
$(M,p,I_A)$ generalises the notion of a Ricci-flat affine manifold.

The model structure is $(\mathbb{S}^n,\si')$ where $\si'$ is the
weight 1 density which arises as a projective polynomial from
$\aI_{\ul{A}}X^{\ul{A}}$ on $\mathbb{R}^{n+1}$ where $\aI_{\ul{A}}$
is a constant covector there. Since $\aI_{\ul{A}}X^{\ul{A}}=0$
describes a hyperplane through the origin in $\mathbb{R}^{n+1}$ it
follows that the zero locus of $\si'$ is a totally geodesic
embedded $\mathbb{S}^{n-1}$ in $\mathbb{S}^{n}$ with its standard
projective structure. The $P$-type decomposition consists of the 3
submanifolds where $\si'$ is, respectively, positive, zero, and
negative. On the open submanifolds $\mathbb{S}^n_{\pm}$ where
$\si'$ is, respectively, positive, and negative, $\si'$ is a
scale and induces the flat connection $\nabla^{\si'}$ in agreement
with the identification (by central projection) of
$\mathbb{S}^n_{\pm}$ with, respectively, the affine subspaces in
$\mathbb{R}^{n+1}$ descibed by $\aI_{\ul{A}}X^{\ul{A}}=\pm 1$.  Thus
these manifolds $(\mathbb{S}^n_{\pm},\nabla^{\si'})$ are
geodesically complete.  According to our general results above these
features are necessarily repoduced in the general situation. Thus we
have the following.
\begin{theorem}\label{std}
  Consider a projective manifold $(M,p)$ equipped with a parallel
  standard cotractor $I_A$.  Then the weight 1 projective density
  $\si=I_AX^A$ satisfies the equation \nn{basic}. The manifold is
  stratified by $P$-types $M_+,M_o,M_-$ according to the
  strict sign of $\si$. These components have a structure as follows:\\
  $\small{\bullet}$ The zero locus $M_o$ of $\si$ is either empty, or
  forms a smooth embedded hypersurface. With respect to any $\nabla\in
  p$, this is totally geodesic, and has canonically an intrinsic
  projective structure $p^{M_o}=[\nabla^{M_o}]$ where $\nabla^{M_o}$
  is simply the restriction of $\nabla$.  The normal tractor
  connection of $(M_o,p^{M_o})$ is naturally a restriction of
  the ambient tractor connection.\\
  $\small{\bullet}$ The open submanifolds $(M_{\pm},\nabla^\si)$ are
  Ricci-flat affine
  manifolds, which are geodesically complete if $M$ is closed.\\
  $\small{\bullet}$ If $M\setminus M_+$ or $M\setminus M_-$ is compact
  (e.g.\ if $M$ is closed this is forced) then $(M\setminus
  M_{\mp},p,I)$ is a geometrically canonical compactification of,
  respectively, $(M_\pm,p)$ (where throughout $p$ and $I$ are
  restricted to the indicated submanifolds).
\end{theorem}
\noindent{\bf Proof:} The last bullet point is simply the observation
that the construction yields a projective analogue of conformal 
compactification of Einstein manifolds.

Concerning the zero locus: On the model the zero locus is a totally
geodesic equatorial embedded $(n-1)$-sphere. It follows at once that
$M_o$ is totally geodesic, as based around a point $q$ in
$\mathcal{Z}(\si)$, by construction $\Phi^{-1}$ is compatible with the
geodesic paths through $q$.  Thus $M_o$ has a projective structure
which is simply a restriction of that from the ambient $(M,p)$.

This result can also be seen via the Thomas cone space since there 
 $\aI_A$ is parallel, and thus is a conormal to the zero locus 
of $\boldsymbol{\si}=\aI_A \zeta^A$ (where we have used
$\and_A\zeta^B=\delta^B_A$ ). Since this zero locus
$\mathcal{Z}(\boldsymbol{\si})$ has a parallel conormal it is totally
geodesic. On the other hand the non-vertical geodesics of $\aM$ are
the lifts of the geodesics from $(M,p)$, 
and it
follows that $M_o$ is totally geodesic. 
Using now that $\mathcal{Z}(\boldsymbol{\si})$ is totally
geodesic, it follows that it inherits an affine manifold structure by
the restriction of the ambient $\and$. The claims about the normal
tractor connection follow, with this restriction of $\and$ to $
\mathcal{Z}(\boldsymbol{\si})$ being the Thomas space over
$(M_o,p|_{M_o})$.

Away from its zero locus, $\si$ is a scale and so we have $\nabla^\si
\si=0$, thus the claim that the components $(M_{\pm},\nabla^\si)$ are
Ricci-flat follows from \nn{basic}.

All other points follow immediately from the corresponding results on
the model via Theorem \ref{canp} (and its proof), Corollary
\ref{zeros}, Theorem \ref{Ptype}, and the discussion of section \ref{prelim}.
\quad $\Box$

\subsection{A parallel tractor metric -- Klein-Einstein structures} 
\label{KleinE}

Here we consider a projective manifold $(M,p)$ equipped with a
non-degenerate symmetric and parallel 2-cotractor $H_{AB}$ of
signature $(r,s)$, $r\geq s\geq 0$. In this case $\si:=H_{AB}X^AX^B$
satisfies the third order equation
\begin{equation}\label{third}
\nabla_{(a}\nabla_{b}\nabla_{c)}\si + 4P_{(ab}\nabla_{c)}\si+ 2\big(\nabla_{(a}P_{bc)}\big)\si=0,
\end{equation}
where $(\cdots)$ indicates the symmeytric part over the enclosed indices.

The model structure is $(\mathbb{S}^n,\si')$ where $\si'$ is the
weight 2 density which arises as a polynomial scalar density
from the homogeneous polynomial
$\boldsymbol{\si}:= \aH_{\ul{A}\ul{B}}X^{\ul{A}}X^{\ul{B}}$ on
$\mathbb{R}^{n+1}$, where $\aH_{\ul{A}\ul{B}}$ is a fixed (signature
$(r,s)$) inner product there.

If $s\geq 1$ then $\aH_{\ul{A}\ul{B}}X^{\ul{A}}X^{\ul{B}}=0$ is a
quadratic variety in $\mathbb{R}^{n+1}$ and, corresponding to this,
the zero locus of $\si'$ is an embebbed variety $S^{r-1}\times S^{s-1} $
in $\mathbb{S}^{n}$ with a signature $(r-1,s-1)$ conformal structure
induced from $\aH_{\ul{A}\ul{B}}$ (viewed now as a metric in
$\mathbb{R}^{n+1}\setminus \{0\}$) restricted to tangent vectors in
$\mathcal{Z}(\boldsymbol{\si})$.  The $P$-type decomposition consists
of the 3 submanifolds where $\si'$ is, respectively, positive,
zero, and negative. On the open submanifolds $\mathbb{S}^n_{\pm}$
where $\si'$ is, respectively, positive, and negative, $\si'$
is a scale and induces a spaceform metric, with signature respectively
$(r,s-1)$ and $(r-1,s)$, and with $\nabla^{\si'}$ the compatible
Levi-Civita connection having curvature, respectively, negative and
positive.  It is well known that these that these manifolds
$(\mathbb{S}^n_{\pm},\nabla^{\si'})$ are geodesically complete.

If $s=0$ then the model is very simple. Then
$\mathcal{Z}(\boldsymbol{\si})$ is empty on $\mathbb{R}^{n+1}\setminus
\{0\}$. There is just the one $P$-type, viz.\ $\mathbb{S}^n$, and, via
$\si$, $H_{AB}$ induces the usual (up to diffeomorphism) unit round metric on this, and
$\nabla^\si$ is the standard Levi-Civita connection.  In general we
have the following.

\begin{theorem}\label{KE}
  Consider a projective manifold $(M,p)$ equipped with a
  non-degenerate symmetric and parallel 2-cotractor $H_{AB}$ of
  signature $(r,s)$, $r\geq s\geq 0$.  The weight 2 projective density
  $\si=H_{AB}X^AX^B$ satisfies the equation \nn{third}.
  If $s=0$ then: \\
  $\small{\bullet}$ $M$ is a single $P$-type and $\si$ is a scale on
  $M$ that induces a positive Einstein metric.\\
  If $s\geq 1$
  we have the following:\\
  $\small{\bullet}$ The manifold is stratified by $P$-types
  $M_+,M_o,M_-$ according to the
  strict sign of $\si$. \\
  $\small{\bullet}$ The zero locus $M_o$ of $\si$ is either, empty, or
  forms a smooth embedded hypersurface with a conformal structure $c$
  of signature $(r-1,s-1)$.  The standard conformal tractor bundle
  agrees with the restriction of the projective tractor bundle $\cT$
  to $M_o$ and the normal conformal tractor connection of $(M_o,c)$ is
  naturally the corresponding restriction of
  the ambient projective tractor connection.\\
  $\small{\bullet}$ On the open submanifold $M_{\pm}$, $\si$ is a
  scale that induces, respectively, a positive/negative Einstein
  metric $g^{\si}$ of signature $(r-1,s)$ or $(r,s-1)$. In each case,
  the affine connection $\nabla^\si$ is the corresponding Levi-Civita
  connection. If
  $M$ is closed then each of $(M_\pm, g^\si)$ is geodesically complete.\\
  $\small{\bullet}$ If $M\setminus M_+$ or $M\setminus M_-$ is compact
  with boundary $M_o$ (e.g.\ if $M$ is closed this is forced) then
  $(M\setminus M_{\mp},p,H)$ is a geometrically canonical
  compactification of, respectively, $(M_\pm, g^\si)$ ($H$ is
  restricted to the indicated submanifolds).
\end{theorem}
\noindent{\bf Proof:} Since a scale trivialises the bundles $\ce(w)$
and splits the dual Euler sequence \nn{ceuler}, it follows that in the
presence of a scale, $H_{AB}$ determines a covariant symmetric two
tensor on the manifold $M$. As in the case of the model, the signature
of this depends on the $P$-type, and indeed via the diffeomorphism
$\Phi$, of Lemma \ref{ncoords3},  we can conclude the signature of
each $P$-type from the model; in either case this is determined in an
obvious way according to whether $X^A$ is timelike or spacelike with
respect to $H_{AB}$. Calculating locally where $\si$ is a scale one sees that
$\nabla^{\cT}_aH_{BC}=0$ implies that $\nabla^\si_a P_{bc}=0$ and that
$P_{bc}$ agrees up to a constant (giving the sign of the scalar
curvature) with the metric induced from $H_{BC}$, and this constant is
non-zero on the open $P$-types \cite{ArmstrongP1,GM}. Since
$\nabla^\si$ is torsion-free it follows that on these $P$-types it is an 
Einstein Levi-Civita connection.

$M_o$ is the set where $X^A$ is null (with respect to $H_{AB}$). But,
where $X^A$ is null, one easily sees that $H_{AB}$ determines a
signature $(r-1,s-1)$ bilinear form, taking values in $\ce(2)|_{M_o}$,
that is independent of any splitting of \nn{ceuler}. This is locally
compatible with the model via $\Phi$. Since $H_{AB}$ is parallel,
this is in particular so along $M_o$, and so $\nabla^\cT$ is metric
preserving along $M_o$.  Here, as elsewhere, we have that in any
choice of weight 1 scale $\si$, and with $\nabla^{\si}$ for the moment
denoting the coupled scale-tractor connection, we have that
$\nabla^{\si}_a(\si^{-1} X^B)$ gives a splitting of
\nn{ceuler}. Combining with the fact that $\nabla^\cT$ is torsion free
it follows, using its characterisation in \cite{CapGoamb}, that
$\nabla^\cT$ agrees with the normal conformal tractor connection along
$M_o$.
 
As in the previous example, all remaining facts follow immediately from 
the corresponding results on
the model via Theorem \ref{canp}, Corollary
\ref{zeros}, Theorem \ref{Ptype}, and the discussion of section \ref{prelim}.
\quad $\Box$

\noindent{\bf Remarks:} The case that $H_{AB}$ has Lorentzian
signature is important. Looking at the model
$(\mathbb{S}^n,\si')$, the part $\mathbb{S}^n_-$ where $\si'$
is negative consists of two copies of hyperbolic space $\mathbb{H}^n$
antipodally placed as the interior of the standard double $S^{n-1}$
quadric on $\mathbb{S}^n$. The compactification $\mathbb{S}^n\setminus
\mathbb{S}^n_+$ adds the boundary spheres. Since this model is based
on central projection (with e.g.\ geodesics arising from planes
through the origin) it is natural to think of the result as two copies
of the Klein model of $\mathbb{H}^n$; whence the curved analogue could
be called a {\em Klein-Einstein manifold} by analogy with the use of
the term Poincar\'e-Einstein in the literature. Note that the
conformal structure of the interior of the Klein-Einstein (KE)
manifold does not extend to the boundary, even though the latter has a
canonical conformal structure. This is clear by continuity
considerations, for example, since the signature of the ambient metric changes
as we cross the zero locus of $\si$. Thus we have a result, which we state as 
 proposition in order to highlight. 

 \begin{proposition} \label{KEP}
A Klein-Einstein manifold involves a
   compactification of its Einstein interior that is strictly
   different to the conformal compactification of a
   Poincar\'e-Einstein (PE) manifold; there is never a smooth
   diffeomporphism between a PE manifold and KE manifold that
   restricts to a conformal map on the interior.
\end{proposition}

Returning to the model with $H_{AB}$ Lorentzian, the component
$\mathbb{S}^n_+$ is the geometry known as de Sitter space in the
general relativity literature; $\mathbb{S}^n\setminus \mathbb{S}^n_-$
is the projective compactification of this. Again $(M\setminus
M_-,p,H_{AB})$ is a curved analogue.

It is not difficult to show that the Klein-Einstein
compactification described here is the same as the {\em projectively
  compact (Einstein) metric} described by Fefferman-Graham in
\cite[chapter 4]{FGnew}. There it is defined as a manifold with
boundary equipped with a negative Einstein metric on the interior that
near the boundary takes form $h/\rho+d{\rho ρ}^2/4\rho^2$; here $\rho$
is a defining function for the boundary, while $h$ a symmetric
2-tensor $h$ which is smooth to the boundary, and with restriction
there a signature $(p,q)$ boundary metric. As explained in
\cite{FGnew}, by an appropriate change of variables these structures
may be transformed to Poincar\'e-Einstein manifolds satisfying an
evenness condition, and hence are also closely related to the
Fefferman-Graham ambient metric. Further details of the geometry of
the Klein-Einstein type structures, and their links to PE manifolds,
is taken up in \cite{Gopro,GM}.

\subsection{Singular and higher codimension zero locus} Examples with
singular zero locus arise easily in the case where we assume more than
one parallel tractor field. For example if we assume $I^1_A$ and
$I^2_A$ are linearly independent parallel cotractor fields, on a given projective manifold $(M,p)$,
then
$S_{AB}:=I^1_AI^2_B+I^1_BI^2_A$ is symmetric and parallel. Thus we are
in the situation of the previous example {\em except} that $S_{AB}$ is
far from non-degenerate. We have
$$
\mathcal{Z}(X^AX^BS_{AB}) = \mathcal{Z}(\si^1\si^2) = \mathcal{Z}(\si^1)\cup \mathcal{Z}(\si^2),
$$
where $\si^1=I^1_AX^A$ and $\si^2=I^2_AX^A$.  In the model, and
generically, this is not smooth.  There are three $P$-types according
to whether none, one, or both of $\si^1$ and $\si^2$ is zero.
Geometrically, on an open dense set, the structure $(M,p,S_{AB}) $ has
projectively related (in the sense of \nn{ptrans}) Ricci-flat affine
structures.

Assume $(M,p)$, $I^1$, and $I^2$ as above, and set
$K_{AB}=I^1_AI^2_B-I^1_BI^2_A$, then we generically obtain a smooth
codimension 2 zero locus 
$$
\mathcal{Z}(k_a)= \mathcal{Z}(\si^1)\cap \mathcal{Z}(\si^2)
$$
for the weight 2 one-form field $k_a = \si^1\nabla_a \si^2 -
\si^2\nabla_a \si^1 $ which corresponds to $K_{AB}X^B$. (Here $\nabla$
is any connection from $p$). We have $K_{AB}X^B=Z_A^b k_b$ where
$Z_A^b$ is the projectively invariant bundle monomorphism in the Euler
sequence \nn{euler}, and note that $\nabla_a \si^{i}$ is non-vanishing
along $\mathbb{Z}(\si^i)$ as $I^i\neq 0$, for $i=1,2$. ) There are two
$P$-types: simply $\mathcal{Z}(k_a)$ and its complement.  Note that
the first BGG equation in this case is
$$
\nabla_{(a}k_{b)}=0.
$$


\begin{thebibliography}{XX}

\bibitem{AndReg} M.T.\ Anderson, {\em Boundary regularity, uniqueness
    and non-uniqueness for AH Einstein metrics on 4-manifolds}, Adv.\
  Math.\ {\bf 179} (2003), 205--249.





\bibitem{A} S.\ Armstrong, {\em Definite signature conformal
    holonomy: a complete classification}, J.\ Geom.\ Phys.\  {\bf 57} (2007),
   2024--2048.

 \bibitem{ArmstrongP1} S.\ Armstrong, {\em Projective holonomy. I.
     Principles and properties}, Ann.\ Global Anal.\ Geom., {\bf 33}
   (2008), 47--69.

 \bibitem{ArmstrongP2} S.\ Armstrong, {\em  Projective
   holonomy. II. Cones and complete classifications},  Ann.\ Global
   Anal.\ Geom.\  {\bf 33} (2008), 137--160.

 \bibitem{BCEG} T.\ Branson, A.\ \v Cap, Eastwood, A.R.\ Gover, {\em
     Prolongations of geometric overdetermined systems}, Internat.\
   J.\ Math.\ {\bf 17} (2006), 641--664.

\bibitem{BEG} T.N.\ Bailey, M.G.\ Eastwood, and A.R.\ Gover, {\em
    Thomas's structure bundle for conformal, projective and related
    structures}, Rocky Mountain J.\ Math.\ {\bf 24} (1994),
  1191--1217.



\bibitem{BGG} I.N.\ Bernstein, I.M.\ Gelfand, S.I.\ Gelfand, {\em
    Differential operators on the base affine space and a study of
    $\frak{g}$-modules}.  Lie groups and their representations
  (Proc.\ Summer School, Bolyai J\'{a}nos Math.\ Soc., Budapest,
  1971), pp. 21--64.

\bibitem{Biquard} O.\ Biquard, M\' etriques d'Einstein asymptotiquement
sym\' etriques,
Ast\' erisque, No. 265 (2000), vi+109 pp.

\bibitem{BiqMazz} O.\ Biquard, R.\ Mazzeo, 
{\em  Parabolic geometries as conformal infinities of Einstein metrics}, 
Arch.\ Math.\ (Brno)  {\bf 42}  (2006),  suppl., 85--104. 

 \bibitem{CQY} A.\ Chang, J.\ Qing, and P.\ Yang,
     {\em Renormalized volumes for conformally compact Einstein manifolds},
     (Russian) Sovrem.\ Mat.\ Fundam.\ Napravl.\  {\bf 17} (2006), 129--142;
     translation in J.\ Math.\ Sci.\ (N.Y.)  {\bf 149} (2008),
     1755--1769

 
\bibitem{CapGoTAMS} A.\ \v Cap, and A.R.\ Gover, {\em Tractor calculi
    for parabolic geometries}, Trans.\ Amer.\ Math.\ Soc.\ {\bf 354}
  (2002), 1511--1548.


 \bibitem{CapGoamb} A.\ \v Cap, and A.R.\ Gover, {\em Standard
     tractors and the conformal ambient metric construction},  Ann.\
     Global Anal.\ Geom.\  {\bf 24} (2003), 231--295.


\bibitem{CGI} A.\ \v Cap, and A.R.\ Gover, {\em CR-Tractors and the
       Fefferman Space}, Indiana University Mathematics Journal, {\bf
       57} (2008), 2519--2570

\bibitem{CGF} A.\ \v Cap, and A.R.\ Gover, {\em A holonomy
       characterisation of Fefferman spaces},  arXiv:math/0611939

\bibitem{CGH} A.\ \v Cap,  A.R.\ Gover, M.\ Hammerl, in progress.

\bibitem{CSbook} A.\ \v Cap, J.\ Slov\'{a}k, Parabolic geometries. I.
  Background and general theory. Mathematical Surveys and Monographs,
  154. American Mathematical Society, Providence, RI, 2009. x+628 pp.

\bibitem{CSSIII} A.\ \v Cap, J.\ Slov\'{a}k, V.\ Sou\v cek,
     {\em Invariant operators on manifolds with almost Hermitian symmetric
       structures. III. Standard operators},  Differential Geom.\ Appl.\
       {\bf 12} (2000), 51--84.

\bibitem{CSS} A.\ \v Cap, J.\ Slov\'{a}k, V.\ Sou\v cek {\em
         Bernstein-Gelfand-Gelfand sequences}, Ann.\ of Math., (2)
       {\bf 154} (2001), 97--113.



\bibitem{CY} S.Y.\ Cheng, S.T.\  Yau, 
{\em On the existence of a complete K\"ahler metric on noncompact 
complex manifolds and the regularity of Fefferman's equation}, 
Comm.\ Pure Appl.\ Math {\bf 33} (1980),  507--544. 


\bibitem{Feff} C. Fefferman, \textit{Monge--Amp\`ere equations,
the Bergman kernel and geometry of pseudoconvex domains},
Ann. of Math. \textbf{103} (1976), 395--416; Erratum \textbf{104} (1976),
393--394. 

\bibitem{FGast} C.\ Fefferman, and C.R.\ Graham, {\em Conformal
    invariants} in: The mathematical heritage of \'{E}lie Cartan (Lyon,
  1984).  Ast\' erisque 1985, Numero Hors Serie, 95--116. 


\bibitem{FGrQ} C.\ Fefferman, and C.R.\ Graham, {\em $Q$-curvature and
    Poincar\'e metrics}, Math.\ Res.\ Lett.\ {\bf 9} (2002), 139--151.

\bibitem{FGnew} C.\ Fefferman, and C.R.\ Graham,  The Ambient Metric,
  Princeton University Press, to appear. \quad arXiv:0710.0919.

\bibitem{FoxIndiana} D.J.F. Fox, {\em Contact projective structures}
  Indiana Univ.\ Math.\ J.\ {\bf 54} (2005), 1547--1598.


\bibitem{FultonH} W.\ Fulton, J.\ Harris, Representation theory. A first course. Graduate Texts in Mathematics, 129. Readings in Mathematics. Springer-Verlag, New York, 1991. xvi+551 pp.

\bibitem{GoPrague} A.R.\ Gover, {\em Almost conformally Einstein
    manifolds and obstructions}, in Differential geometry and its
  applications, 247--260, Matfyzpress, Prague, 2005. Electronic:
  arXiv:math/0412393

\bibitem{GoMGC} A.R.\ Gover, {\em Conformal Dirichlet-Neumann maps and
    Poincar\'e-Einstein Manifolds}, SIGMA (Symmetry, Integrability and
  Geometry: Methods and Applications)  {\bf 3} 100 (2007).
\quad arXiv:0710.2585


\bibitem{Go-al} A.R. Gover, {\em Almost Einstein and
    Poincar\'e-Einstein manifolds in Riemannian signature},   J.\
    Geometry and Physics, {\bf 60} (2010), 182--204.

\bibitem{Gopro}  A.R.\ Gover, in progress.

\bibitem{GoLeitprog} A.R.\ Gover, and F.\ Leitner, {\em A class of
    compact Poincare-Einstein manifolds: properties and construction},
  {\em Communications in Contemporary Mathematics}, to appear. \quad
  arXiv:0808.2097



\bibitem{GM}  A.R.\ Gover, and H. Macbeth, in progress. 

\bibitem{GrSrni} C.R.\ Graham, {\em Volume and area renormalizations
    for conformally compact Einstein metrics}. The Proceedings of the
    19th Winter School "Geometry and Physics" (Srn\'\i, 1999), Rend.\
    Circ.\ Mat.\ Palermo (2) Suppl. No. 63 (2000), 31--42.



 \bibitem{GrL} C.R.\ Graham, and J.M.\ Lee, {\em Einstein metrics with
    prescribed conformal infinity on the ball}, Adv.\ Math.\ {\bf 87} (1991),
     186--225. 

   \bibitem{GrZ} C.R.\ Graham, and M.\ Zworski, {\em Scattering matrix
       in conformal geometry}, Invent.\ Math.\ {\bf 152} (2003), 89--118.

\bibitem{GrW} C.R.\ Graham, and E.\ Witten, {\em Conformal anomaly of
    submanifold observables in AdS/CFT correspondence}, Nuclear Phys.
  B 546 (1999), no. 1-2, 52--64.

 \bibitem{GKP} S.S.\ Gubser, I.R.\ Klebanov and A.\ M.\ Polyakov,
  {\em Gauge theory correlators from non-critical string theory},
  Phys.\ Letters B {\bf 428} (1998), 105--114.

\bibitem{KMS} I.\ Kolar, P.W.\ Michor, J.\ Slov\'ak, {\em Natural
    Operations in Differential Geometry}, Springer 1993

\bibitem{lee} J.M.\ Lee, {\em Fredholm operators and Einstein metrics
on conformally compact manifolds}, Mem.\ Amer.\ Math.\ Soc.\ {\bf 183}
(2006), no. 864, vi+83 pp.


\bibitem{L} F.\ Leitner, {\em  On transversally
    symmetric pseudo-Einstein and Fefferman-Einstein spaces}, Math.\
  Z.\ {\bf 256} (2007), 443--459.

\bibitem{Lepowsky} J.\ Lepowsky, {\em A generalization of the
    Bernstein-Gelfand-Gelfand resolution}, J.\ Algebra, {\bf 49}
  (1977), 496--511.

\bibitem{Maldacena} J.\ Maldacena {\em The large $N$ limit of superconformal
    field theories and supergravity},  Adv.\ Theor.\ Math.\ Phys.\  {\bf 2}
    (1998), 231--252.


\bibitem{MP} R.\ Mazzeo and F.\ Pacard, {\em Maskit combinations 
of Poincar\'e-Einstein metrics},  Adv.\ Math.\  {\bf 204}  (2006),  379--412.

\bibitem{Orsted} B.\ Orsted, {\em Generalized gradients and Poisson
    transforms}, Global analysis and harmonic analysis
  (Marseille-Luminy, 1999), 235--249, S\'{e}min. Congr., 4, Soc. Math.
  France, Paris, 2000.

\bibitem{ot} R. Penrose and W. Rindler,
{\em Spinors and Space-time vol.~1},
Cambridge University Press 1984.

\bibitem{Rthesis} M.\ Randall, Almost projectively Ricci-flat
  manifolds, MSc thesis, University of Auckland, 2010.

\bibitem{SS} K.\ Skenderis, and S.N.\ Solodukin, {\em Quantum effective
action from the AdS/CFT correspondence}, Phys.\ Lett.\ {\bf B472}
(2000), 316--322, hep-th/9910023

\bibitem{Thomas} T.Y. Thomas, {\em Announcement of a projective theory of affinely connected manifolds}, Proc.\ Nat.\ Acad.\ Sci., {\bf 11} (1925), 588--589.





\end{thebibliography}
\end{document}